\documentclass[11pt,a4paper]{article}
\usepackage[T1]{fontenc}
\usepackage{lmodern,amsmath,amssymb,amsthm,mathtools,microtype}
\usepackage[margin=25mm]{geometry}
\usepackage[shortlabels]{enumitem}
\usepackage{booktabs,array,needspace,tabularx}
\usepackage{tikz}
\usetikzlibrary{arrows.meta,calc,positioning}
\usepackage[hidelinks]{hyperref}
\usepackage{fancyhdr}
\setlist[enumerate]{topsep=4pt,itemsep=2pt}
\setlist[itemize]{topsep=4pt,itemsep=2pt}
\newtheorem{theorem}{Theorem}[section]
\newtheorem{maintheorem}{Theorem}

\newtheorem{proposition}[theorem]{Proposition}
\newtheorem{lemma}[theorem]{Lemma}
\newtheorem{corollary}[theorem]{Corollary}
\theoremstyle{definition}
\newtheorem{definition}[theorem]{Definition}
\newtheorem{remark}[theorem]{Remark}
\newcommand{\F}{\mathbb F}

\newcommand{\Prob}{\mathbb P}

\newcommand{\supp}{\operatorname{supp}}

\hypersetup{bookmarksdepth=2,pdftitle={Hyperbolicity and zero divisors in random graphical groups},pdfauthor={Elena I. Bunina, Winston Lee, Alexei Miasnikov, and Kyrylo Muliarchyk}}
\title{Hyperbolicity and zero divisors in random graphical groups}
\author{Elena I. Bunina, Winston Lee, Alexei Miasnikov, and Kyrylo Muliarchyk}
\date{}
\newcommand{\inv}[1]{\overline{#1}}
\newcommand{\Z}{\mathbb{Z}}
\newcommand{\R}{\mathbb{R}}
\newcommand{\one}{\mathbf{1}}
\newcommand{\cP}{\mathcal{P}}
\newcommand{\cA}{\mathcal{A}}
\newcommand{\cC}{\mathcal{C}}
\newcommand{\cD}{\mathcal{D}}

\newcommand{\cR}{\mathcal{R}}

\newcommand{\cU}{\mathcal{U}}
\newcommand{\chain}{\mathsf{c}}
\newcommand{\dmin}{d_{\min}}
\newcommand{\dmax}{d_{\max}}
\newcommand{\pmin}{p_{\min}}
\newcommand{\wmin}{w_{\min}}
\newcommand{\Dcl}{D_{\mathrm{cl}}}
\newcommand{\deff}{d_{\mathrm{eff}}}
\newcommand{\norm}[1]{\lVert #1\rVert}
\newcommand{\Dnorm}{\norm{D}}
\newcommand{\bd}{\partial}
\newcommand{\tP}{\mathbb{P}}
\DeclareMathOperator{\rk}{rk}
\DeclareMathOperator{\Bij}{Bij}
\DeclareMathOperator{\girth}{girth}

\DeclareMathOperator{\dist}{dist}
\DeclareMathOperator{\St}{St}
\newtheorem{example}[theorem]{Example}
\numberwithin{equation}{section}
\begin{document}
\maketitle
\begin{abstract}
A labelled graph defines a group on its alphabet by taking the labels of closed paths as relators. We prescribe the proportions of the possible sets of outgoing labels, match occurrences of each generator uniformly to occurrences of its inverse, and condition on girth at least $\lfloor c\log n\rfloor$, where $n$ is the number of vertices and $c>0$ is sufficiently small. For a balanced distribution of outgoing labels, a matrix records the weights of pairs of paths that read the same word. If its spectral radius is less than one and each vertex has degree at least three, then, with probability tending to one, the group is non-elementary hyperbolic and has a finite two-dimensional classifying space. On the same event, every path from a prescribed root to a different vertex has nontrivial label in the group. We prove the linear isoperimetric inequality by subdividing diagrams and applying a local-to-global theorem. For OpenAI's characteristic-two zero-divisor construction, the criterion recovers the hyperbolicity result of Mghirbi and Douglas.
\end{abstract}

\begingroup
\setlength{\parskip}{0pt}
\setcounter{tocdepth}{1}
\small
\tableofcontents
\endgroup
\newpage
\section{Introduction}\label{sec:intro}

Let $k$ be a field. Kaplansky asked whether $k[G]$ can contain two nonzero elements whose product is zero when $G$ is torsion-free \cite{kaplansky57}. Torsion gives such a product immediately: if $g$ has order $s>1$, then $(1-g)(1+g+\cdots+g^{s-1})=0$. For two nonzero factors to have zero product, the contributions to every coefficient must cancel. An order on $G$ invariant under multiplication on both sides rules this out: the product of the largest elements of the two supports has a unique representation, so its coefficient is nonzero.

Labelled graphs provide a way to arrange repeated products. Rips and Segev used them to construct torsion-free groups without the unique product property \cite{rips-segev}; Steenbock proved hyperbolicity for such constructions by graphical small cancellation \cite{steenbock}. Repetition alone does not force a coefficient to vanish. In OpenAI's construction over $\F_2$, the outgoing labels make the multiplicities even, and the geometry keeps the two factors nonzero \cite{seed}. OpenAI's subsequent construction gives a failure of direct finiteness: there are elements $a,b$ with $ab=1$ and $ba\ne1$ \cite{dfseed}. Odd-characteristic constructions appear in \cite{ivanov2026,blmm2026,md2026,jaikin2026}. For a different torsion-free hyperbolic group, Douglas and Mghirbi prove non-soficity by ruling out approximate permutation actions \cite{dmnonsofic}.

We study the geometric question suggested by these examples. Which distributions of outgoing labels guarantee, with high probability, that the group is hyperbolic and that every path from a prescribed root to a different vertex has nontrivial label? The criterion below and its proof for general outgoing-label distributions are due to Winston Lee. The algebraic cancellation in the characteristic-two application is due to OpenAI. Mghirbi and Douglas already proved hyperbolicity for that application \cite[Corollary~5.3]{md2026}.

\subsection{Scope and main theorem}\label{intro:scope}

Let $S$ be a set of $r\ge2$ generators and write $T=S\sqcup S^{-1}$, with inversion $a\mapsto\bar a$. For a vertex of the graph, its \emph{link} is the set of labels on its outgoing edges. A probability distribution $\pi$ with rational values prescribes the proportion of vertices with each link $\tau\subseteq T$. For letters $a,b\in T$, put
\[
 p_a=\pi\{\tau:a\in\tau\},\qquad
 p_{a,b}=\pi\{\tau:a,b\in\tau\}\quad(a\ne b).
\]
Assume that the distribution is balanced and uses every letter:
\[
\text{(H1)}\qquad p_a=p_{\bar a}>0\quad\text{for every }a.
\]
The turn matrix $W_\pi=(w_\pi(a,b))$ and the collision matrix $M_\pi$ are defined by
\[
 w_\pi(a,b)=
 \begin{cases}p_{\bar a,b}/p_{\bar a},& b\ne\bar a,\\[1pt]0,& b=\bar a,\end{cases}
 \qquad M_\pi=W_\pi\circ W_\pi,
\]
where $\circ$ denotes entrywise multiplication. For $b\ne\bar a$, the entry $w_\pi(a,b)$ is the fraction of links containing $\bar a$ that also contain $b$. The remaining hypotheses are
\[
\begin{array}{ll}
\text{(H2)}&|\tau|\ge3\quad\text{whenever }\pi(\tau)>0,\\
\text{(H3)}&\rho(M_\pi)<1,
\end{array}
\]
where $\rho$ is the spectral radius.
Let $d_{\max}=\max\{|\tau|:\pi(\tau)>0\}$,
$p_{\min}=\min_a p_a$, and
\[
 c_*(\pi)=\min\left\{
 \frac1{\log(2|T|/p_{\min})},\frac1{2\log d_{\max}}\right\}.
\]
Call $n$ admissible if every $n\pi(\tau)$ is an integer. Assign exactly $n\pi(\tau)$ vertices the link $\tau$. For each $s\in S$, choose a uniform bijection from the vertices with outlabel $s$ to those with outlabel $s^{-1}$, independently for different $s$. Write $\Gamma_\Phi$ for the resulting labelled graph. Fix $0<c<c_*(\pi)$ and condition on $\Gamma_\Phi$ having girth at least $L=\lfloor c\log n\rfloor$. Throughout the graphical model, \emph{with high probability} means that the conditional probability tends to one as $n$ tends to infinity through admissible values.

The group $G_\Phi$ has generators $S$ and relators the labels of all closed paths in $\Gamma_\Phi$. Choose a spanning tree in each component. The labels of the associated fundamental cycles form a finite set of defining relators $B_\Phi$. Let $X_\Phi$ be the rose on $S$ with a cone attached along the label map of each component.

\Needspace{14\baselineskip}
\begin{maintheorem}[Outlabel criterion]\label{intro:construction}
Assume \textup{(H1)--(H3)}. For each admissible $n$, fix a root $x_0$ before sampling the graph. There is a constant $\varepsilon_1=\varepsilon_1(\pi,c)>0$ such that the following statements hold simultaneously with high probability:
\begin{enumerate}
\item every word $v$ over $T$ that represents the identity in $G_\Phi$ has area at most $|v|/\varepsilon_1$ over $\langle S\mid B_\Phi\rangle$;
\item no nonempty cyclically reduced word of length at most
      $\varepsilon_1L$ is trivial in $G_\Phi$;
\item $X_\Phi$ is a finite two-dimensional classifying space for $G_\Phi$, and
      $G_\Phi$ is torsion-free and non-elementary hyperbolic;
\item every path from $x_0$ to a different vertex has nontrivial label
      in $G_\Phi$.
\end{enumerate}
The same conclusions hold for finitely many independent parts with
fixed positive proportions and separate matching laws, as specified
in Section~\ref{sec:partitioned}.
\end{maintheorem}

All three hypotheses can be checked from $\pi$ before sampling; none requires a parity condition. Assertion~(4) separates the chosen root from the other vertices of its component. It does not assert that paths from the root to two different non-root vertices have distinct labels in the group. For OpenAI's two link distributions, root separation keeps the two factors nonzero, while the parity identity gives their zero product in $\F_2[G_\Phi]$ \cite{seed}.

\subsection{Plan of proof}\label{intro:plan}

A word can label two paths only if the required turns are available on both. For words of a given length, the sum of their squared weights is the sum of the entries of a power of the collision matrix. Hypothesis (H3) makes this sum decay exponentially (Section~\ref{sec:criterion}). Switching two matched edges lets us use these estimates after conditioning on girth. The graph estimates also bound the diameter of each component (Section~\ref{sec:estimates}).

Next express a diagram in terms of \emph{stars}: the triangles at one cone vertex form a face whose perimeter reads a path in the graph (Section~\ref{sec:diagrams}). A diagram with little exposed boundary forces many label coincidences. We sort the repeated intervals by how often they occur; the squared-weight estimate makes the probability of such a bounded pattern tend to zero (Section~\ref{sec:local}).

The last step removes the bound on diagram size. Short chords cut each long star into tiles. The reassembly argument lets us apply the local estimate to the subdisks required by the local-to-global theorem (Section~\ref{sec:global}). The resulting inequality bounds the area of a disk by a constant times its boundary length and excludes reduced spherical diagrams with faces. A diagram with one distinguished root face rules out a trivial label from the root to a different vertex (Section~\ref{sec:proofs}). Sections~\ref{sec:partitioned} and~\ref{sec:graph-application} apply these conclusions to independent parts and to the zero-divisor construction. Section~\ref{sec:comparison} gives two restrictions on zero divisors in the different model where relator words are sampled independently.

\section{The model}\label{sec:model}

\subsection{Alphabet, words and graphs}

Fix a finite set $S$ with $r:=|S|\ge2$, and put $T:=S\sqcup S^{-1}$. The fixed-point-free involution $a\mapsto\inv a$ exchanges each $s$ with $s^{-1}$. A word $a_1\cdots a_h$ over $T$ is \emph{reduced} if $a_{i+1}\ne\inv{a_i}$ for all $i<h$. A nonempty reduced word is \emph{cyclically reduced} if $a_1\ne\inv{a_h}$. Let $F(S)$ be the free group on $S$ and $R_S$ the rose with one edge for each $s\in S$.

Graphs are finite and may have loops and multiple edges. Each edge has two orientations, denoted $e$ and $\inv e$. A \emph{$T$-labelled graph} assigns a letter $\ell(e)\in T$ to each oriented edge, with $\ell(\inv e)=\inv{\ell(e)}$. Its label map to $R_S$ is an \emph{immersion} if the oriented edges leaving any vertex have distinct labels. A path $e_1\cdots e_h$ has length $h$ and label $\ell(e_1)\cdots\ell(e_h)$. It is \emph{non-backtracking} if $e_{i+1}\ne\inv{e_i}$ for all $i<h$; a nonempty closed non-backtracking path is \emph{cyclically non-backtracking} if also $e_1\ne\inv{e_h}$. Loops and pairs of parallel edges are cycles of length $1$ and $2$, respectively. The \emph{girth} is the least length of a cycle, and is $\infty$ for a forest. The \emph{cycle rank} is $\rk:=|E|-|V|+\#\{\text{components}\}$. Distances are path distances, with distance $\infty$ between different components.

\begin{lemma}\label{lem:closed-nb}
A closed non-backtracking path of positive length contains a cycle; in particular its length is at least the girth.
\end{lemma}

\begin{proof}
Suppose the image of the path is a forest. This image is connected and has an edge, so it is a tree with at least two vertices. Measure distances within this tree from the initial vertex, and choose a vertex $z$ at maximal distance $d\ge1$. Each neighbour of $z$ has distance $d-1$ or $d+1$, and exactly one has distance $d-1$. Maximality excludes distance $d+1$, so $z$ has degree one in the image. Since $z$ is not the initial vertex, the path enters and leaves $z$ consecutively along the same edge, a backtrack. This contradiction shows that the image contains a cycle. The path traverses every edge of its image, so its length is at least the length of that cycle, and hence at least the girth.
\end{proof}

The \emph{link} of a vertex is the set of oriented edge-germs leaving it (equivalently the sphere of radius $1/4$ in the metric graph whose edges have length one). A loop contributes two germs. In an immersed labelled graph the labels of these germs are distinct, so we identify the link with its set of outgoing labels, or \emph{outlabels}.

\subsection{Links and the turn matrix}

\begin{definition}[Link distributions]\label{def:outlabel}
An \emph{outlabel set}, also called a \emph{link}, is a subset $\tau\subseteq T$. A \emph{link distribution} is a probability measure $\pi$ on the subsets of $T$ taking rational values. Put $p_a:=\pi\{\tau:a\in\tau\}$ and, for $a\ne b$, $p_{a,b}:=\pi\{\tau:a,b\in\tau\}$. Let $\dmin$ and $\dmax$ be the smallest and largest cardinalities of links in the support of $\pi$, and let $N_\pi$ be the least positive integer such that $N_\pi\pi(\tau)\in\Z$ for every $\tau$. Call $\pi$ \emph{balanced} if $p_a=p_{\inv a}$ for all $a$. Say that $\pi$ \emph{uses every letter} if $p_a>0$ for all $a$; in this case put $\pmin:=\min_ap_a$.
\end{definition}

\begin{definition}[Turn matrix, collision matrix, word weights]\label{def:turn}
Let $\pi$ be balanced and use every letter. The \emph{turn matrix} $W_\pi=(w_\pi(a,b))_{a,b\in T}$ is
\[
w_\pi(a,b):=\frac{p_{\inv a,b}}{p_{\inv a}}\quad(b\ne\inv a),\qquad w_\pi(a,\inv a):=0,
\]
and the \emph{collision matrix} is its entrywise square $M_\pi:=W_\pi\circ W_\pi$. The \emph{weight} of a word $v=a_1\cdots a_h$ with $h\ge1$ is $P(v):=\prod_{i=1}^{h-1}w_\pi(a_i,a_{i+1})$, with the empty product equal to one when $h=1$. When $W_\pi$ has a positive entry, as it does under (H2) below, let $\wmin$ be its least positive entry and put $a_0:=\log(1/\wmin)$. We write $w$ for $w_\pi$.
\end{definition}

If a path reaches a vertex $x$ by reading $a$, then $\inv a$ is an outlabel of $x$. For $b\ne\inv a$, the entry $w(a,b)$ is the conditional $\pi$-probability that a link contains $b$, given that it contains $\inv a$. It therefore measures how frequently a non-backtracking continuation by $b$ is available among the prescribed links. This interpretation does not assume that successive turns along a sampled path are independent.

\begin{lemma}\label{lem:identities}
Let $\pi$ be balanced and use every letter.
\begin{enumerate}[(1)]
\item $\sum_{v\in T^h}P(v)^2=\one^{\mathsf T}M_\pi^{\,h-1}\one$ for $h\ge1$, where $\one$ is the all-ones vector.
\item $0\le w\le1$, $w(a,\inv a)=0$, and $p_a\,w(a,b)=p_b\,w(\inv b,\inv a)$ for all $a,b$.
\item $\sum_bw(a,b)=\mathbb E_\pi\bigl[\,|\tau|-1\bigm|\inv a\in\tau\,\bigr]\ge\dmin-1$.
\item If $P(v)>0$, then $v$ is reduced.
\end{enumerate}
\end{lemma}

\begin{proof}
(1) Expanding the matrix product gives one summand for each word in $T^h$. (2) The bounds follow from the definition. If $b\ne\inv a$, both sides of the reversal identity equal $p_{\inv a,b}$ by balance; if $b=\inv a$, both vanish. (3) Divide $\sum_{b\ne\inv a}p_{\inv a,b}=\sum_{\tau:\inv a\in\tau}\pi(\tau)(|\tau|-1)$ by $p_{\inv a}$. (4) An immediate cancellation contributes the factor $w(a,\inv a)=0$.
\end{proof}

\subsection{The random graph}

\begin{definition}[Configurations]\label{def:graph}
Let $\pi$ be balanced and use every letter, and let $n$ be \emph{admissible}, that is, a positive multiple of $N_\pi$. Put $V_n:=\{1,\dots,n\}$ and fix an \emph{outlabel assignment} $\tau_n$ from $V_n$ to the subsets of $T$ with $|\tau_n^{-1}(\tau)|=n\pi(\tau)$ for every $\tau$. For $a\in T$ let $V_a:=\{x\in V_n:a\in\tau_n(x)\}$; thus $|V_a|=np_a$, and $|V_s|=|V_{s^{-1}}|$ by balance. The \emph{configuration space} is
\[
\Omega_n:=\prod_{s\in S}\Bij\bigl(V_s,V_{s^{-1}}\bigr).
\]
For $\Phi=(\phi_s)_{s\in S}\in\Omega_n$, the graph $\Gamma_\Phi$ has vertex set $V_n$ and an edge $(s,x)$ for each $s\in S$ and $x\in V_s$, joining $x$ to $\phi_s(x)$ and labelled $s$ in this direction.
\end{definition}

Only $\Phi$ is random, and all estimates depend on $\tau_n$ only through the numbers $n\pi(\tau)$.

\begin{lemma}[Immersion]\label{lem:immersion}
For $x\in V_n$ and $a\in T$, exactly one oriented edge of $\Gamma_\Phi$ leaves $x$ with label $a$ if $a\in\tau_n(x)$, and none otherwise. Hence the label map $\Gamma_\Phi\to R_S$ is an immersion, and $\deg x=|\tau_n(x)|$, a loop counting twice.
\end{lemma}

\begin{proof}
An edge leaving $x$ with label $s\in S$ is $(s,x)$, which exists if and only if $x\in V_s$. An edge leaving $x$ with label $s^{-1}$ is $(s,\phi_s^{-1}(x))$ read backwards, which exists if and only if $x\in V_{s^{-1}}$.
\end{proof}

A \emph{prescription} is a triple $(x,s,y)$ with $s\in S$, $x\in V_s$ and $y\in V_{s^{-1}}$; it \emph{holds} in $\Phi$ if $\phi_s(x)=y$. The edges of $\Gamma_\Phi$ correspond bijectively to the prescriptions that hold.

\subsection{Girth conditioning}

For a balanced distribution $\pi$ using every letter and satisfying $\dmax\ge2$, put
\[
c_*(\pi):=\min\Bigl\{\frac1{\log(2|T|/\pmin)},\ \frac1{2\log\dmax}\Bigr\}.
\]
Fix $0<c<c_*(\pi)$, put $L=L_n:=\lfloor c\log n\rfloor$ and $\Omega_n^c:=\{\Phi\in\Omega_n:\girth(\Gamma_\Phi)\ge L\}$, and let $\mu_n$ be the uniform probability measure on $\Omega_n^c$. Under (H1) and (H2) below, $\Omega_n^c\ne\emptyset$ for all large admissible $n$ (Lemma~\ref{lem:girth-exists}). An event holds \emph{with high probability} (w.h.p.) if its $\mu_n$-probability tends to $1$ as $n\to\infty$ through admissible values. This is the conditioning of~\cite{seed,dfseed}, with $c_*(\pi)$ in place of~\cite[(2.5)]{seed}.

\subsection{The group}

\begin{definition}[The group and the cone complex]\label{def:group}
Let $G_\Phi:=\langle S\mid\ell(\gamma):\gamma\text{ a closed path in }\Gamma_\Phi\rangle$. The \emph{cone complex} $X_\Phi$ is the rose $R_S$ together with a cone on each component of $\Gamma_\Phi$, attached along the label map.
\end{definition}
Choose a spanning tree and a base vertex in each component. Let $p_x$ be the tree path from the base vertex of the component of $x$ to $x$. An oriented edge $e$ from $x$ to $x'$ outside the trees defines the fundamental cycle $c_e=p_xe\inv{p_{x'}}$. Let $B_\Phi$ be the finite set of their labels. Proposition~\ref{prop:cone} shows that $\pi_1(X_\Phi)\cong G_\Phi=\langle S\mid B_\Phi\rangle$. The elements of $B_\Phi\cup B_\Phi^{-1}$ are called the relators of this presentation.

\begin{remark}\label{rem:parameter}
The model depends on $\pi$ and $c$, but hypothesis (H3) below depends only on $W_\pi$. By Lemma~\ref{lem:identities}, a turn matrix is nonnegative, vanishes at the entries $(a,\inv a)$, satisfies the reversal identity $p_aw(a,b)=p_bw(\inv b,\inv a)$ for a positive vector $p$ with $p_a=p_{\inv a}$, and has row sums at least $2$ under (H2); we do not characterize the matrices that arise. The OpenAI construction~\cite{seed} runs two independent copies of the model, with two link distributions forcing parity; each copy alone is an instance (Example~\ref{ex:openai}).
\end{remark}

\section{The collision criterion}\label{sec:criterion}

The hypotheses on the link distribution are:
\begin{itemize}
\item[\textbf{(H1)}] $\pi$ is balanced and uses every letter;
\item[\textbf{(H2)}] $\dmin\ge3$;
\item[\textbf{(H3)}] $\rho(M_\pi)<1$.
\end{itemize}
Balance ensures that each letter has as many prescribed occurrences as its inverse, so the required bijections exist. The assumption that every letter is used excludes generators that do not occur on graph edges. The graph need not be connected: the diameter estimate applies to each component, and the cone complex has one cone per component.

The conclusions of Theorem~\ref{intro:construction} follow from these hypotheses.
The last assertion concerns a root fixed before sampling; it does not assert simultaneous injectivity of the entire graph in the Cayley graph. The constant $\varepsilon_1$ can be chosen as in~\eqref{eq:epsone}; it depends only on $\pi$ and $c$, not on $n$.

\begin{lemma}[Forms of (H3)]\label{lem:H3forms}
For a nonnegative $T\times T$ matrix $M$ the following are equivalent:
\begin{enumerate}[(1)]
\item $\rho(M)<1$;
\item there are $\delta>0$ and $h_0$ with $\one^{\mathsf T}M^{h-1}\one\le e^{-2\delta h}$ for all $h\ge h_0$;
\item there are $f:T\to(0,\infty)$ and $\lambda<1$ with $Mf\le\lambda f$ entrywise.
\end{enumerate}
For $M=M_\pi$, the quantity in \textup{(2)} is $\sum_{|v|=h}P(v)^2$ by Lemma~\ref{lem:identities}.
\end{lemma}

\begin{proof}
Let $\norm{\cdot}$ be the maximum row sum norm. Since $M$ is nonnegative,
\[
\norm{M^k}\le\one^{\mathsf T}M^k\one\le|T|\norm{M^k}.
\]
Gelfand's formula gives $\rho(M)=\lim_{k\to\infty}\norm{M^k}^{1/k}$~\cite[Cor.~5.6.14]{HJ13}, so (2) implies $\rho(M)\le e^{-2\delta}<1$. Conversely, assume (1) and choose $\rho(M)<\lambda<1$. The same formula gives $\norm{M^{h-1}}\le\lambda^{h-1}$ for all sufficiently large $h$. Hence
\[
\one^{\mathsf T}M^{h-1}\one\le|T|\lambda^{h-1}\le\lambda^{h/2}
\]
for large $h$, which is (2) with $\delta:=-\frac14\log\lambda$.

Under (1), the vector $f:=\sum_{k\ge0}M^k\one$ is well defined and has positive entries. It satisfies $Mf=f-\one\le(1-1/\max f)f$, proving (3). Finally, (3) implies $M^kf\le\lambda^kf$ and therefore $\norm{M^k}\le\lambda^k\max f/\min f$. Taking $k$th roots gives $\rho(M)\le\lambda<1$.
\end{proof}

\begin{corollary}\label{cor:rowtest}
\textup{(H3)} holds if
\[
\max_a\sum_bw(a,b)^2<1.
\]
In particular, it holds whenever
\[
\Bigl(\max_bw(a,b)\Bigr)\Bigl(\sum_bw(a,b)\Bigr)<1
\qquad\text{for every }a\in T.
\]
\end{corollary}

\begin{proof}
The spectral radius is at most the maximum row sum, so
\[
\rho(M_\pi)\le\norm{M_\pi}=\max_a\sum_bw(a,b)^2.
\]
For each row,
\[
\sum_bw(a,b)^2\le\Bigl(\max_bw(a,b)\Bigr)\Bigl(\sum_bw(a,b)\Bigr).
\]
\end{proof}

By Lemma~\ref{lem:identities}(3), $\sum_bw(a,b)+1$ is the mean degree of the vertices with outlabel $\inv a$.

\begin{definition}[Effective density]\label{def:deff}
Put $\deff(\pi):=\frac12\bigl(1+\log\rho(M_\pi)/\log(2r-1)\bigr)$, with $\log0:=-\infty$. Thus (H3) holds if and only if $\deff(\pi)<\frac12$.
\end{definition}

\begin{example}[Uniform outlabels]\label{ex:uniform}
Let $3\le D\le2r$ and let $\pi$ be uniform on the $D$-element subsets of $T$. Every letter has probability $D/(2r)$ of belonging to the link, so (H1) holds; (H2) follows from $D\ge3$. Conditional on one letter belonging to the link, each of the other $2r-1$ letters belongs with probability $(D-1)/(2r-1)$. Hence $w(a,b)=\kappa:=(D-1)/(2r-1)$ for $b\ne\inv a$. Thus $M_\pi=\kappa^2(J-\Pi)$, where $J$ is the all-ones matrix and $\Pi$ is the permutation matrix of $a\mapsto\inv a$. The nonnegative matrix $J-\Pi$ has every row sum equal to $2r-1$, so
\[
\rho(M_\pi)=\frac{(D-1)^2}{2r-1},\qquad\deff(\pi)=\frac{\log(D-1)}{\log(2r-1)}.
\]
Theorem~\ref{intro:construction} therefore applies precisely when $(D-1)^2<2r-1$. For example, $(r,D)=(3,3)$ gives $\rho=\frac45$ and $(8,4)$ gives $\frac35$, while $(2,3)$ gives $\frac43$ and $(5,4)$ gives $1$.
The first two examples have vertex degrees three and four, respectively, below the minimum degree five in the fixed-graph criterion of Mghirbi and Douglas \cite[Definition~3.3]{md2026}.
\end{example}

In Gromov's density model, random groups at density $d<\frac12$ are w.h.p.\ infinite and hyperbolic, and at $d>\frac12$ they are w.h.p.\ trivial or of order $2$~\cite{ollivier,invitation}. The normalization above places the collision condition at the same numerical threshold $1/2$. It does not establish a phase transition at $\deff=1/2$ in the outlabel model.

\begin{example}[The OpenAI outlabels]\label{ex:openai}
Let $q=128$, $v=q^2+q+1$ and $p=(q+1)/v$. In~\cite[\S2]{seed}, $T$ consists of the $v$ points of the projective plane over $\F_q$ and three extra letters, each paired by the involution with an ordinary point. For $\pi_B$, choose a projective line uniformly and, independently, adjoin all three extra letters with probability $p$, and none otherwise. For $\pi_A$, choose a line uniformly. Independently of the line, choose one of the three two-letter subsets of the extra letters, each with probability $p/2$, or the empty subset with probability $1-3p/2$, and adjoin it to the line.

Every letter occurs with probability $p$ under either distribution, so (H1) holds. Both distributions have $\dmin=q+1$, giving (H2). For $\pi_B$ and a non-cancelling turn $b\ne\inv a$, the entry $w(a,b)$ equals $1/(q+1)$ if $\inv a$ and $b$ are both ordinary, $p$ if exactly one is extra, and $1$ if both are extra. Let $f(a)=5$ if $\inv a$ is extra and $f(a)=1$ otherwise. If $\inv a$ is ordinary, at most three of the $v-1=q(q+1)$ ordinary successors $b$ have $f(b)=5$, and the three extra successors have $f(b)=1$. If $\inv a$ is extra, then $f(a)=5$ and the two other extra successors contribute $2$. Hence
\[
(M_\pi f)(a)\le\frac q{q+1}+\frac{12}{(q+1)^2}+3p^2\quad\text{or}\quad(M_\pi f)(a)\le2+(v+12)p^2,
\]
respectively, and both bounds are less than $\frac{149}{150}f(a)$; this is the estimate of~\cite[Lemma~2.1]{seed}. Thus (H3) holds by Lemma~\ref{lem:H3forms}(3). For $\pi_A$, the turn probabilities are the same except that the value $1$ for two distinct extra letters becomes $\frac12$, so the same vector $f$ works. Theorem~\ref{intro:construction} therefore applies to each distribution separately. Their joint presentation is covered by Theorem~\ref{thm:partitioned}.
\end{example}

\section{Estimates for the random graph}\label{sec:estimates}

Throughout Sections~\ref{sec:estimates}--\ref{sec:proofs}, we assume (H1) and (H2) and fix $0<c<c_*(\pi)$. Hypothesis (H3) enters the probabilistic argument in Section~\ref{sec:local}; the later geometric conclusions use that local estimate. Constants depend only on $\pi$ and $c$ unless stated otherwise. ``For large $n$'' means for all admissible $n\ge n_0$, where $n_0$ depends only on $\pi$, $c$ and the named parameters. The notation $o(L)$ denotes a quantity of absolute value at most $\zeta_nL$, where $\zeta_n\to0$ depends only on the same data. For large $n$,
\[
\frac c2\log n\le L\le c\log n,\qquad L\ge3.
\]
Put $r_n:=8\dmax^{2L}$. Since $2c\log\dmax<1$, we have $r_n=o(n)$.

\subsection{Switching}

Let $\tP$ be the uniform probability measure on $\Omega_n$, under which the $\phi_s$ are independent and uniform.

\begin{lemma}\label{lem:uniform-presc}
If $k<n\pmin$, then any $k$ distinct prescriptions hold simultaneously with $\tP$-probability at most $(n\pmin-k)^{-k}$.
\end{lemma}

\begin{proof}
If the prescriptions are incompatible, their probability is zero. Otherwise, the $j_s$ prescriptions concerning $\phi_s$ fix its values at $j_s$ distinct points, with probability $(N_s-j_s)!/N_s!\le(n\pmin-k)^{-j_s}$, where $N_s:=|V_s|\ge n\pmin$. Independence and $\sum_sj_s=k$ give the bound after multiplying over $s$.
\end{proof}

\begin{lemma}[Switching]\label{lem:switch}
Let $\Phi\in\Omega_n$ and $s\in S$, and let $(x,s,y)$ and $(u,s,z)$ hold in $\Phi$ with $\dist(\{x,y\},\{u,z\})>2L$ in $\Gamma_\Phi$. Let $\Phi'$ agree with $\Phi$ except that $\phi'_s(x)=z$ and $\phi'_s(u)=y$. Then $\Phi'\in\Omega_n$, the graphs $\Gamma_{\Phi'}$ and $\Gamma_\Phi$ differ only in the edges $(s,x)$ and $(s,u)$, and if a cycle of $\Gamma_{\Phi'}$ of length less than $L$ contains one of the \emph{new} edges $(s,x)$, $(s,u)$ of $\Gamma_{\Phi'}$, then the edge $(s,u)$ of $\Gamma_\Phi$ lies on a cycle of $\Gamma_\Phi$ of length less than $L$.
\end{lemma}

\begin{proof}
The distance condition gives $x\ne u$, so $\phi'_s$ is a bijection. Let $\Gamma_0$ be $\Gamma_\Phi$ without the edges $(s,x)$ and $(s,u)$; it is also $\Gamma_{\Phi'}$ without its new edges. Removing the new edges from a cycle $C$ of $\Gamma_{\Phi'}$ of length less than $L$ leaves one or two paths in $\Gamma_0$, of total length less than $L$, joining in pairs the endpoints $x,z$ and $u,y$ of the new edges in $C$. A path from $\{x,y\}$ to $\{u,z\}$ would contradict the distance condition. Hence $C$ contains both new edges, and one of the paths joins $z$ to $u$. Being part of $C$, it is simple, and it avoids the edge $(s,u)$ of $\Gamma_\Phi$; with that edge it forms a cycle of $\Gamma_\Phi$ of length less than $L$.
\end{proof}

\subsection{Large girth}

\begin{lemma}\label{lem:girth-exists}
Put $\beta:=c\log(2|T|/\pmin)$, so that $0<\beta<1$. For large $n$, $\tP(\Omega_n^c)\ge\exp(-n^{(1+\beta)/2})$; in particular $\Omega_n^c\ne\emptyset$.
\end{lemma}

\begin{proof}
Call an edge \emph{short} if it lies on a cycle of length less than $L$, and let $X(\Phi)$ be the number of short edges of $\Gamma_\Phi$. A cycle of length $k$ is traced by a sequence $(x_0,a_1,x_1,\dots,a_k,x_k=x_0)$ of vertices and letters, and there are at most $(2r)^kn^k$ such sequences; its $k$ edges are distinct prescriptions. By Lemma~\ref{lem:uniform-presc}, for large $n$,
\[
\mathbb E_{\tP}X\le\sum_{k<L}k(2r)^kn^k(n\pmin-k)^{-k}\le L^2\Bigl(\frac{2|T|}{\pmin}\Bigr)^L\le L^2n^\beta.
\]
Since $X$ is integer-valued, Markov's inequality gives $\tP(X\le m_0)\ge\frac12$ for $m_0:=\lfloor2L^2n^\beta\rfloor$. Also $m_0<n\pmin/4$ for large $n$.

Let $0<X(\Phi)\le m_0$, and let $(x,s,y)$ be the prescription of a short edge. At most $4\dmax^{2L}$ vertices lie within distance $2L$ of $\{x,y\}$. Each is an endpoint of at most two edges labelled $s$, so all but at most $r_n$ of the $np_s\ge n\pmin$ edges labelled $s$ are at distance more than $2L$ from $\{x,y\}$. Since fewer than $n\pmin/4$ edges are short and $r_n=o(n)$, for large $n$ we can choose such a distant edge $(u,s,z)$ that is not short. By Lemma~\ref{lem:switch}, no cycle of length less than $L$ in the switched graph $\Gamma_{\Phi'}$ contains a new edge. Every other cycle is already a cycle of $\Gamma_\Phi$. Thus no new short edge is created, and the short edge $(s,x)$ is removed.

Repeat the procedure, using a fixed rule for the choices. The number of short edges decreases at every step, so we reach some $\Phi^*\in\Omega_n^c$ after $t\le X(\Phi)\le m_0$ switchings. A switching is undone from its triple $(s,x,u)$, since $\phi_s(x)=\phi'_s(u)$ and $\phi_s(u)=\phi'_s(x)$, and there are at most $|S|n^2$ possible triples. Hence $\Phi$ is determined by $\Phi^*$ and a sequence of $t\le m_0$ triples, giving
\[
\tfrac12|\Omega_n|\le\bigl|\{\Phi:X(\Phi)\le m_0\}\bigr|\le(m_0+1)\,(|S|n^2)^{m_0}\,|\Omega_n^c|.
\]
As $\log\bigl(2(m_0+1)(|S|n^2)^{m_0}\bigr)=O(n^\beta\log^3n)$, this gives the claim for large $n$.
\end{proof}

\begin{lemma}[Conditioned prescriptions]\label{lem:cond-presc}
Let $\mathcal S$ be a set of $k$ prescriptions that hold with positive $\mu_n$-probability, and let $e=(x,\sigma,y)\notin\mathcal S$ be a prescription with $np_\sigma>k+r_n$. Then
\[
\mu_n(e\text{ holds}\mid\mathcal S\text{ holds})\le\frac1{np_\sigma-k-r_n}.
\]
Consequently:
\begin{enumerate}[(a)]
\item if $m+r_n\le n\pmin/2$, any $m$ distinct prescriptions hold simultaneously with probability at most $(2/(n\pmin))^m$;
\item for each $C'>0$ and large $n$, distinct prescriptions $e_i=(x_i,\sigma_i,y_i)$, $1\le i\le E\le C'L$, hold simultaneously with probability at most $e^{o(L)}\prod_i(np_{\sigma_i})^{-1}$, where $o(L)$ depends only on $\pi,c,C'$.
\end{enumerate}
\end{lemma}

\begin{proof}
Let $\Omega_{\mathcal S}\subseteq\Omega_n^c$ be the set of configurations in which $\mathcal S$ holds, and $\Omega_{\mathcal S,e}\subseteq\Omega_{\mathcal S}$ those in which $e$ holds too. Let $\Phi\in\Omega_{\mathcal S,e}$. As in the previous proof, at least $np_\sigma-k-r_n$ prescriptions $(u,\sigma,z)\notin\mathcal S$ hold in $\Phi$ and satisfy $\dist(\{x,y\},\{u,z\})>2L$. For each of them the switch $\Phi'$ of Lemma~\ref{lem:switch} has girth at least $L$, because $\Gamma_\Phi$ does, and satisfies $\mathcal S$; so $\Phi'\in\Omega_{\mathcal S}$. Since $z=\phi'_\sigma(x)$ and $u=\phi'^{-1}_\sigma(y)$, the pair $(\Phi,(u,\sigma,z))$ is determined by $\Phi'$. Hence $|\Omega_{\mathcal S,e}|\,(np_\sigma-k-r_n)\le|\Omega_{\mathcal S}|$, which is the claim.

For (a) and (b), the conclusion is immediate if the joint event has probability zero. Otherwise all preceding conditioning events have positive probability, and the chain rule writes the joint probability as $\prod_i\mu_n(e_i\text{ holds}\mid e_1,\dots,e_{i-1}\text{ hold})$. In (a), each factor is at most $(n\pmin-(i-1)-r_n)^{-1}\le2/(n\pmin)$. In (b), each factor is at most $(np_{\sigma_i})^{-1}(1-x_n)^{-1}$, where $x_n:=(C'L+r_n)/(n\pmin)\to0$. Since $(1-x_n)^{-E}\le e^{2x_nC'L}=e^{o(L)}$ once $x_n\le\frac12$, the claimed estimate follows.
\end{proof}

\subsection{Expansion and diameter}
For $U\subseteq V_n$, let $B_j(U)$ be the vertices at distance at most $j$ from $U$. The following estimate gives exponential growth of small balls and a logarithmic bound on the diameter of each component.
\begin{lemma}[Expansion and diameter]\label{lem:expansion}
There are $C_0\ge1$ and $d_0>0$, depending only on $\pi,c$, such that the following hold simultaneously with high probability:
\begin{enumerate}[(a)]
\item $|B_1(U)|\ge\frac54|U|$ whenever $1\le|U|\le n/C_0$;
\item every component of $\Gamma_\Phi$ has diameter at most $d_0L$.
\end{enumerate}
\end{lemma}
\begin{proof}
Choose $C_0$ large, depending only on $\pi$, so that $2/C_0<\pmin/2$ and the constants in the estimates below satisfy $O(1/C_0)\le1$ and $O(C_0^{-1/4})\le1/2$. Suppose $|U|=k\le n/C_0$ and $|B_1(U)|<5k/4$. Extend $B_1(U)$ to a set $U'$ of size $y=\lceil5k/4\rceil$. Since every vertex has degree at least three, at least $q=\lceil3k/2\rceil$ edges have an endpoint in $U$, and all these edges lie in $U'$. If $y<L$, the induced graph on $U'$ is a forest, so it has at most $y-1<q$ edges, a contradiction. Thus $k>4(L-1)/5$, and in particular $k\ge48$ for large $n$.
Count the choices of $U'$ and of $q$ prescriptions in it. Lemma~\ref{lem:cond-presc}(a) gives
\[
\begin{aligned}
\mu_n(\text{failure of size }k)
&\le\binom ny\binom{ry^2}q\left(\frac2{n\pmin}\right)^q\\
&\le\left(\frac{en}k\right)^{5k/4+1}\left(O(k/n)\right)^{3k/2}\\
&=\frac{en}k\left(O((k/n)^{1/4})\right)^k.
\end{aligned}
\]
Here $k\le y\le2k$, $q\ge3k/2$, and the base $O(k/n)$ is at most one. If $k\le\sqrt n$, this is at most $en\,n^{-k/16}\le en^{-2}$ for large $n$. If $\sqrt n<k\le n/C_0$, it is at most $e\sqrt n\,2^{-\sqrt n}$. Summing over at most $n$ values proves (a).
For (b), put $R=\lceil\log n/\log(5/4)\rceil$. Repeated application of (a) shows that every ball of radius $R$ has more than $n/C_0$ vertices. Along a geodesic of length $d$, the radius-$R$ balls centred at positions $0,2R+1,2(2R+1),\ldots$ are disjoint. There are more than $d/(2R+1)$ such balls, and each has more than $n/C_0$ vertices. Hence $d<(2R+1)C_0\le d_0L$, where $d_0=(4/(c\log(5/4))+3)C_0$ works for large $n$.
\end{proof}

\subsection{Short closures}

\begin{lemma}[Short closures]\label{lem:closure}
Let $\Gamma$ be a $T$-labelled graph with at most $n$ vertices in each component, immersed in $R_S$, without loops or multiple edges, of minimum degree at least $3$ and component diameters at most $d_0L$, where $L\ge\max\{1,\frac c2\log n\}$. Put $\Dcl:=d_0+8/(c\log2)+12$. Then every non-backtracking path $\sigma$ of positive length extends to a closed, cyclically non-backtracking path $\sigma\rho$ with $1\le|\rho|\le\Dcl L$.
\end{lemma}

\begin{proof}
\emph{Based loops.} Let $f$ be an unoriented edge incident to a vertex $z$, and put $k:=\lceil\log_2n\rceil+2$. Since every vertex has degree at least three, at least $2^{k-1}>n$ non-backtracking walks of length $k$ start at $z$ along $f$. Two distinct such walks $P,Q$ therefore end at the same vertex. Remove their longest common terminal segment. The remaining paths $P''$ and $Q''$ are nonempty, since the original walks are distinct and have equal length. They still start along $f$, and their last edges differ. Thus $J_z(f):=P''\inv{Q''}$ is a non-backtracking closed path at $z$, of length at most $2k$, that leaves and returns along $f$.

\emph{Closure.} Let $\sigma$ run from $a$ to $b$, with first edge $e_1$ and last edge $e_\ell$. If $a\ne b$, choose a geodesic $R$ from $b$ to $a$, with first edge $r_1$ and last edge $r_2$. Choose an incident unoriented edge $f_b$ at $b$ different from the underlying edges of $e_\ell$ and $r_1$, and similarly choose $f_a$ at $a$ different from the underlying edges of $r_2$ and $e_1$. Minimum degree three ensures that these choices exist. Put $\rho:=J_b(f_b)\,R\,J_a(f_a)$. If $a=b$, choose an incident edge $f_a$ different from the underlying edges of $e_\ell$ and $e_1$, and put $\rho:=J_a(f_a)$.

Each constituent path is non-backtracking, and at every junction, including the cyclic junction, the incoming and outgoing underlying edges are distinct. Thus $\sigma\rho$ is cyclically non-backtracking. The added path is nonempty, and
\[
|\rho|\le d_0L+4k\le\Dcl L,
\]
because $k\le\frac2{c\log2}L+3$ and $L\ge1$.
\end{proof}

\section{Diagrams over a triangular presentation}\label{sec:diagrams}

The diagram reductions and bounded-pattern method originate in
\cite[Sections~3--5]{seed}. For the original designs, the boundary
accounting that gives hyperbolicity is also proved by
Mghirbi--Douglas \cite[Lemma~4.3 and Proposition~4.4]{md2026}.
Here these arguments are organized using links, triangular diagrams,
and stars, so that the collision condition can be used for an
arbitrary prescribed link distribution.

In this section, $\Gamma=\Gamma_\Phi$ for a fixed $\Phi$ with $\girth(\Gamma)\ge3$, so $\Gamma$ has no loops or multiple edges. We also fix a root $x_0\in V_n$, chosen independently of $\Phi$. We express the diagram arguments of~\cite{seed} in the terminology of Ollivier~\cite{ollivier}.

\subsection{Diagrams}

Let $Z:=\{z_x:x\in V_n\}$ be a set of new letters. For an oriented edge $e$ of $\Gamma$ from $x$ to $x'$ with label $a$, put $r_e:=z_x^{-1}\,a\,z_{x'}$, so that $r_{\inv e}=r_e^{-1}$. Fix an orientation $e^+$ of each edge and let $\cP_\Phi:=\langle S\cup Z\mid r_{e^+}\ (e\in E(\Gamma))\rangle$. Since $x\ne x'$, the three letters of $r_e$ are distinct and no two are mutually inverse. For $y\ne x_0$, call $z_y^{-1}z_{x_0}$ a \emph{root relator}, and let $\cP_\Phi^+$ be $\cP_\Phi$ with these root relators added. They will be used only in the proof of root separation.

\begin{definition}[Diagrams~{\cite[\S1.3]{ollivier}}]\label{def:diagram}
Let $\langle\cA\mid\cR\rangle$ be a presentation.
\begin{enumerate}[(i)]
\item A \emph{(possibly filamentous) diagram} is a finite connected planar combinatorial $2$-complex $D\subset\R^2$ whose oriented edges carry letters of $\cA^{\pm}$, opposite orientations carrying inverse letters, and each of whose $2$-cells (\emph{faces}) $c$ carries a relator $r_c\in\cR$, a marked vertex, and an orientation, such that $\bd c$, read from the marked vertex in that orientation, spells $r_c$. Only the interior of a face is required to be an open disk.
\item The \emph{holes} of $D$ are the bounded components of $\R^2\setminus D$, and $D$ is \emph{$0$-hole} if it has none. The \emph{boundary words} of $D$ are read along the frontiers of the components of $\R^2\setminus D$, an edge with the complement on both sides being read twice; $|\bd D|$ is their total length.
\item A \emph{spherical diagram} carries the same decorations on a cell structure of $S^2$; then $|\bd D|:=0$.
\item Two distinct faces form a \emph{reducible pair} if they share an edge, bear the same relator with opposite orientations, and the common edge has the same position in that relator, counted from their marked vertices. $D$ is \emph{reduced} if it has no reducible pair. $|D|$ denotes the number of faces.
\end{enumerate}
\end{definition}

We use two standard facts~\cite[Ch.~V]{LS77}: (vK) a word is a product of at most $N$ conjugates of relators and their inverses in the free group if and only if it is the boundary word of a $0$-hole, possibly filamentous, diagram with at most $N$ faces; in particular it represents $1$ in $\langle\cA\mid\cR\rangle$ if and only if such a diagram exists; and (Min) among such diagrams with a given boundary word, one with the fewest faces is reduced. We fix an orientation of $\R^2$ (or of $S^2$), read face boundaries \emph{positively}, that is, with the face on the left, and call a face \emph{positive} if its orientation is the positive one.

\subsection{Stars}
\begin{definition}[Root faces and $T$-diagrams]\label{def:Tdiagram}
A \emph{root face} bears a root relator. A \emph{$T$-diagram} is a reduced planar or spherical diagram over $\cP_\Phi^+$, with any number of holes in the planar case, having at most one root face and boundary words over $T$.
\end{definition}
In this subsection $D$ is a $T$-diagram. Reverse the ambient orientation if necessary so that its root face, if present, is positive. A $Z$-edge is oriented towards its \emph{head} when it reads some $z_x$; its other end is its \emph{tail}. The remaining edges are $T$-edges. Sectors between consecutive edge-ends at a vertex that lie in faces are called corners. An \emph{apex} is a marked vertex of a face, and its \emph{star} $\St(o)$ is the union of $o$, the open faces marked at $o$, and the open $Z$-edges at $o$.
A triangle $c$ bearing $r_{e^+}$ lies over $e_c=e^+$ if positive and over $\inv{e^+}$ otherwise. Its positive boundary reading from its apex is $z_x^{-1}az_{x'}$, where $e_c$ goes from $x$ to $x'$ with label $a$. The corresponding reading of a root face is $z_y^{-1}z_{x_0}$. Thus
\begin{equation}\label{eq:corners}
\begin{minipage}{0.88\textwidth}
a corner contains a head if and only if it is the marked corner of its face, and then both edge-ends are heads. Every other corner consists of a $Z$-tail and a $T$-edge-end, or is the unmarked corner of the root face and consists of two tails.
\end{minipage}
\end{equation}
\begin{lemma}[Stars]\label{lem:stars}
\begin{enumerate}[(a)]
\item Every sector at an apex is a marked corner, and every edge-end there is a head.
\item The triangles $c_1,\ldots,c_k$ at an apex $o$, read in positive cyclic order starting just after the root face when it is marked at $o$, are distinct. Their edges $e_{c_1},\ldots,e_{c_k}$ form a non-backtracking path $\gamma_o$ in $\Gamma$. If the root face is marked at $o$ and bears $z_y^{-1}z_{x_0}$, then $\gamma_o$ goes from $x_0$ to $y$ and $k\ge1$. Otherwise it is closed and cyclically non-backtracking, and $k\ge\girth(\Gamma)$.
\item The $T$-edges of these triangles form a closed path $\bd_o$. The star $\St(o)$ is open, and is the image of the interior of a closed disk $Q_o$ under a characteristic map $q_o:Q_o\to D$ which is a homeomorphism on the interior and maps $\bd Q_o$ along $\bd_o$.
\item Faces on opposite sides of a $T$-edge are distinct triangles. The edges of $\Gamma$ over which they lie have inverse labels and are distinct even up to orientation.
\end{enumerate}
\end{lemma}
\begin{proof}
(a) Start at a marked corner. A sector adjacent to a head must be a corner, since no $Z$-edge occurs in a boundary word. By~\eqref{eq:corners} it is marked and its other end is a head. Continuing around $o$ proves the assertion.
(b) Every corner at $o$ is marked, so each incident face occurs once. The last edge of each face and the first edge of the next are the same $Z$-edge, traversed in opposite directions in their positive readings. Away from the root face, their labels therefore make the triangles' edges concatenate in $\Gamma$. If the root face is marked at $o$, its two distinct letters force at least one triangle and identify the path's endpoints as $x_0,y$. Otherwise the path closes. A consecutive reversal would give triangles with the same relator in opposite orientations and the common $Z$-edge in the same relator position, a reducible pair. Thus there is no backtrack, including at the cyclic join in the non-root case. The girth gives the stated length bound.
(c) No $Z$-edge has both endpoints at $o$, because every edge-end there is a head. Glue the incident faces cyclically around a vertex. This gives a disk with boundary the $T$-edges of the triangles, including in the root case, where the intervening digon closes the fan. The characteristic maps give the required map. Its restriction to the interior is injective, since the open faces and edges are disjoint and embedded; it is a homeomorphism onto the open star.
(d) The root face has no $T$-edge, and each triangle uses its unique $T$-edge once. The two positive readings traverse the shared edge in opposite directions. If the corresponding graph edges agreed up to orientation, the triangles would form a reducible pair along their middle letter.
\end{proof}
\begin{definition}[Facets and traces]\label{def:facets}
The \emph{facets} of a star $f=\St(o)$ are the occurrences of edges in $\bd_o$; write $|\bd f|=k$. Order them linearly after choosing a first triangle; for the root star, start just after the root face. The \emph{trace} $\gamma_f=\gamma_o$ sends the facet of $c_i$ to $e_{c_i}$. Two facets on the same $T$-edge are \emph{paired}; the others are \emph{free}. Put
\[
\Dnorm=\sum_f|\bd f|,\qquad \operatorname{free}(D)=\#\{\text{free facets}\}.
\]
Thus $\Dnorm$ counts triangles, and is $|D|$ or $|D|-1$ according as there is no root face or one. Distinct free facets lie on distinct boundary edges, so $\operatorname{free}(D)\le|\bd D|$, and $\operatorname{free}(D)=0$ for spherical diagrams.
\end{definition}

\begin{definition}[Arcs and segments]\label{def:arcs}
Let $D_T$ be the graph formed by the $T$-edges lying on at least one face, together with their endpoints. Following Gruber~\cite[\S2.2]{Gru15}, an \emph{arc} is a path in $D_T$, embedded except possibly at its endpoints, whose interior vertices have degree $2$ in $D_T$ and whose endpoints, which may coincide, do not; it is \emph{interior} if each of its edges lies on two faces, and $\cA(D)$ denotes the set of interior arcs. A vertex of $D_T$ is a \emph{breakpoint} if its degree in $D_T$ is not $2$ or the first facet of some star begins there. The \emph{segments} are the paths into which the breakpoints cut $D_T$; a segment is \emph{interior} if each of its edges lies on two faces. An \emph{interval} of a star is a nonempty set of consecutive facets in its linear order.
\end{definition}

\begin{lemma}[Segments]\label{lem:segments}
\begin{enumerate}[(a)]
\item Let $v$ be a vertex of degree $2$ in $D_T$, with edge-ends $a,a'$ in $D_T$. On each of the two sides of $v$ between $a$ and $a'$, either the edges of $a$ and $a'$ both have an incident face on that side, or neither does. In the first case, their facets on that side are consecutive facets of one star, unless the first facet of that star begins at $v$.
\item Every component of $D_T$ contains a breakpoint. On each side of an interior segment with $j$ edges, the facets along it form an interval of $j$ facets of one star, and the $i$-th facet of one of these intervals is paired with the $(j+1-i)$-th facet of the other. Every paired facet lies along exactly one interior segment.
\item The number of interior segments is at most $|\cA(D)|$ plus the number of stars.
\end{enumerate}
\end{lemma}

\begin{proof}
(a) By Lemma~\ref{lem:stars}(a), $v$ is not an apex; so by~\eqref{eq:corners} every corner at $v$ consists of a $Z$-tail and an end of an edge of $D_T$, except possibly the unmarked corner of the root face, which has two $Z$-tails. Both sectors beside a $Z$-edge-end are corners, since $Z$-letters occur in no boundary word, while the $T$-edges at $v$ outside $D_T$ have no face on either side. Turn around $v$ from $a$ towards $a'$. If the first sector is a corner, its other edge-end is a $Z$-tail, and the following corners, possibly including the root face, lead to the other $T$-edge-end $a'$. If no root corner intervenes, the two corners belong to faces sharing a $Z$-edge, which by Lemma~\ref{lem:stars}(b) are consecutive faces $c_i,c_{i+1}$ of one star, and their facets are consecutive unless $c_{i+1}$ is the first face. If the unmarked corner of the root face intervenes, there are instead three corners: those of the last triangle, the root face and the first triangle of the root star. In that case its first facet begins at $v$, which is precisely the allowed exception. If the first sector is not a corner, no $Z$-end is adjacent to it, so up to $a'$ only filament ends occur and no sector on that side is a corner. Then neither edge lies on a face on that side, since such a face would have a corner at $v$ there.

(b) Each edge of $D_T$ lies on some closed path $\bd_o$, which passes through the breakpoint where the first facet of its star begins; so every component contains a breakpoint and every segment ends at breakpoints. At the interior vertices of an interior segment no first facet begins, so by (a) the facets on each side form an interval of one star. The two stars traverse the segment in opposite directions, which gives the reversed pairing. A paired facet lies on a $T$-edge on two faces; this edge lies on exactly one segment, which is interior by (a).

(c) By (a), each interior segment lies in an interior arc or in a component of $D_T$ that is a cycle of vertices of degree $2$, all of whose edges lie on two faces. The interior segments arise by cutting these arcs and cycles at the vertices where first facets begin, and each such cycle contains one by (b). Cutting a path at $j$ interior vertices gives $j+1$ paths, and a cycle at $j\ge1$ vertices gives $j$; as there is one first facet per star, the bound follows.
\end{proof}

\subsection{Patterns}

\begin{definition}[Patterns and realizations]\label{def:pattern}
Let $D$ be as above, with a first face chosen in each star. The \emph{trace graph} $\Delta\subseteq\Gamma$ is the union of the traces $\gamma_f$. For a star $f$ let $\lambda_f$ be $\gamma_f$ regarded as a path in $\Delta$, closed unless $f$ is the root star, and for a facet $\xi$ of $f$ put $\lambda(\xi):=\gamma_f(\xi)$. If a root face is present, mark the initial vertex of its trace as the root of $\Delta$. The \emph{pattern} of $D$ is $\cD:=(D,\Delta,\lambda)$. Two patterns, possibly over different configurations, are \emph{equivalent} if there are a bijection $f\mapsto f'$ of their stars and an isomorphism $\iota$ of their unlabelled trace graphs with $\iota\circ\lambda_f=\lambda'_{f'}$, preserving the root star when present and carrying the pairs of intervals along interior segments to each other (each pair recorded by the two stars, the positions of the first facets of the two intervals, and the length). A \emph{realization} of $\cD$ in $\Gamma_\Phi$ is a graph map $\psi:\Delta\to\Gamma_\Phi$, injective on vertices and edges and commuting with $e\mapsto\inv e$, mapping the root to $x_0$ when present and such that $\psi(\lambda(\xi))$ and $\psi(\lambda(\xi'))$ carry mutually inverse letters whenever $\xi$ and $\xi'$ are paired.
\end{definition}

\begin{lemma}\label{lem:pattern-props}
\begin{enumerate}[(a)]
\item Each $\lambda_f$ is a non-backtracking path, closed and cyclically non-backtracking unless it is the root trace, of length $|\bd f|$ in $\Delta$; $\Delta$ is the union of their images, and $\girth(\Delta)\ge\girth(\Gamma)$.
\item If $\xi$ and $\xi'$ are paired, then $\lambda(\xi)$ and $\lambda(\xi')$ carry mutually inverse letters and are distinct edges of $\Delta$, even up to orientation.
\item The inclusion $\Delta\subseteq\Gamma$ realizes $\cD$. The numbers of stars and facets, the pairing of facets, $\operatorname{free}(D)$, the number of interior segments, and the existence of a realization in a given $\Gamma_{\Phi'}$ depend only on the equivalence class of $\cD$.
\end{enumerate}
\end{lemma}

\begin{proof}
(a) and (b) restate Lemma~\ref{lem:stars}(b),(d). In (c), the pairing is determined by the pairs of intervals along interior segments (Lemma~\ref{lem:segments}(b)), and the other quantities by the pairing and the paths $\lambda_f$.
\end{proof}

\section{A local isoperimetric inequality}\label{sec:local}

In this section, (H1)--(H3) hold, and $\delta>0$ and $h_0$ are as in Lemma~\ref{lem:H3forms}(2) for $M=M_\pi$, so that $\sum_{|v|=h}P(v)^2\le e^{-2\delta h}$ for $h\ge h_0$. Put
\begin{equation}\label{eq:eps}
\varepsilon:=\min\Bigl\{\frac14,\ \frac{\delta}{16(1+a_0)}\Bigr\}.
\end{equation}

\begin{theorem}[Local isoperimetric inequality]\label{thm:local}
Let $K\ge1$ and $I\ge0$ be integers and $C\ge1$. With high probability, every $T$-diagram $D$ with at most $K$ stars, at most $I$ interior arcs, and $L\le\Dnorm\le CL$ satisfies $\operatorname{free}(D)>\varepsilon\Dnorm$; in particular $|\bd D|>\varepsilon\Dnorm$.
\end{theorem}

We fix $K$, $I$, $C$ until the end of the section; $o(L)$ and ``large $n$'' may depend on them. We work on the event $\girth(\Gamma_\Phi)\ge L$, which has $\mu_n$-probability one. The proof follows Ollivier's scheme~\cite[\S2.2]{ollivier} as adapted in~\cite[\S3]{seed}.

\subsection{Trace graphs are thin}

\begin{lemma}\label{lem:thin}
Let $D$ be as in Theorem~\ref{thm:local}, with pattern $(D,\Delta,\lambda)$.
\begin{enumerate}[(a)]
\item $\rk\Delta\le c_0:=C+r_0$, where $r_0:=\max\{r\ge1:r\le4C(\lfloor\log_2r\rfloor+2)\}$.
\item Every vertex of $\Delta$ has degree at least $1$, and only the two endpoints of the root trace can have degree $1$. Let $V_\star$ be the set of vertices of $\Delta$ that have degree other than $2$ or are endpoints of some $\lambda_f$, and put $c_1:=6c_0+2K+10$. Then $\sum_{z\in V_\star}\deg z\le c_1$, so $|V_\star|\le c_1$, and every component of $\Delta$ meets $V_\star$. Suppressing the vertices outside $V_\star$ gives a multigraph $\Theta$ on $V_\star$ with at most $c_1/2$ edges, each corresponding to a non-backtracking path in $\Delta$ (a \emph{chain}) with interior vertices outside $V_\star$; the chains partition the edges of $\Delta$.
\item Each $\lambda_f$ is a concatenation of \emph{traversals} of chains from end to end. The multiplicity $m(e):=\#\{\xi:\lambda(\xi)\in\{e,\inv e\}\}$ is constant along chains and at most $2(C+K)$, and there are at most $c_2:=c_1(C+K)$ traversals.
\end{enumerate}
\end{lemma}

\begin{proof}
(a) Since $\Delta$ is the union of the traces, it has at most $\Dnorm\le CL$ edges, and its girth is at least $L$. Repeatedly deleting vertices of degree at most $1$ preserves the cycle rank and leaves a graph $\Delta^\circ$ of minimum degree at least $2$. Its components that are cycles have at least $L$ edges each, so there are at most $C$ of them. Suppressing the vertices of degree $2$ in the other components gives a multigraph $\Xi$ of minimum degree at least $3$, with $r$ edges, say, whose lengths (those of the corresponding paths) sum to at most $CL$, and with at most $2r/3$ vertices. Also, $\rk\Delta\le C+r$.

We may assume $r\ge1$. Consider the non-backtracking random walk $(\tau_i)$ on the $2r$ oriented edges of $\Xi$, started uniformly: $\tau_{i+1}$ is uniform among the $\deg z-1\ge2$ oriented edges issuing from the endpoint $z$ of $\tau_i$ other than $\inv{\tau_i}$. The uniform distribution is stationary, since each oriented edge issuing from $z$ has $\deg z-1$ predecessors, each moving to it with probability $1/(\deg z-1)$. So $\tau_1,\dots,\tau_k$ have expected total length at most $kCL/r$, and with probability at least $\frac12$ the total length is at most $2kCL/r$. Each walk has probability at most $(2r)^{-1}2^{1-k}$, so there are at least $r2^{k-1}$ walks of total length at most $2kCL/r$. For $k:=\lfloor\log_2r\rfloor+2$ this exceeds $r^2$, which exceeds the number of ordered pairs of vertices of $\Xi$; so two distinct such walks have the same endpoints. Expanding them gives distinct non-backtracking paths $P\ne Q$ in $\Delta$ with the same endpoints, since distinct oriented edges of $\Xi$ at a vertex begin with distinct edges of $\Delta$. Removing their longest common terminal segment leaves $P'$ and $Q'$, and $P'\inv{Q'}$ is a closed non-backtracking path of positive length (if both $P',Q'$ are nonempty, their last edges differ). By Lemma~\ref{lem:closed-nb}, $L\le|P|+|Q|\le4kCL/r$, so $r\le4C(\lfloor\log_2r\rfloor+2)$ and $r\le r_0$.

(b) Every trace is non-backtracking, and all but the root trace are closed and cyclically non-backtracking. Hence every vertex of $\Delta$ has degree at least two except possibly the two endpoints of the root trace. Since $\sum_z(\deg z-2)\le2c_0$, the sum of $\deg z-2$ over vertices of degree at least three is at most $2c_0+2$. Their total degree is at most $6c_0+6$. Add at most two vertices of degree one and at most $K+1$ marked endpoints of degree two: the degree sum over $V_\star$ is at most $c_1=6c_0+2K+10$. Each component contains an initial vertex of a trace. Thus each edge-end at $V_\star$ begins a chain, and the chains partition the edges; there are at most $c_1/2$ chains.

(c) A non-backtracking path entering a chain runs through it to the other end, and each $\lambda_f$ starts and ends in $V_\star$; this gives the traversals and the constancy of $m$ along chains. If $\lambda_f$ traverses an oriented edge at positions $i<i'$, its subpath from position $i$ to $i'-1$ is a closed non-backtracking path, so $i'-i\ge L$ by Lemma~\ref{lem:closed-nb}. Hence $\lambda_f$ traverses each oriented edge at most $|\bd f|/L+1$ times, and $m(e)\le2(\Dnorm/L+K)\le2(C+K)$. The number of traversals is $\sum_{\text{chains}}m\le\frac{c_1}2\cdot2(C+K)=c_2$.
\end{proof}

\begin{lemma}\label{lem:count}
There is $c_3$, depending only on $K,I,C$, such that for large $n$ the patterns of diagrams as in Theorem~\ref{thm:local} fall into at most $L^{c_3}$ equivalence classes.
\end{lemma}

\begin{proof}
By Lemmas~\ref{lem:thin} and~\ref{lem:segments}, the class is determined by: the multigraph $\Theta$, the root star if present, and, for each star, the sequence of chain traversals of $\lambda_f$, for which there are boundedly many possibilities; the lengths of the at most $c_1/2$ chains, each at most $CL$; and, for each of the at most $I+K$ interior segments, its two stars, the positions of the first facets of its two intervals, and its length, for which there are at most $K^2(CL)^3$ possibilities. This gives at most $O(1)\,(CL)^{c_1/2}(1+K^2(CL)^3)^{I+K}\le L^{c_3}$ classes.
\end{proof}

\subsection{Stages}

For $j\ge1$ let $\Delta_j$ consist of the edges $e$ of $\Delta$ with $m(e)\ge j$ and their endpoints, with $V_j$ vertices and $E_j$ edges. Thus $\Delta_1=\Delta$, and $\Delta_j=\emptyset$ for $j>m_\star:=\max_e m(e)$, where $m_\star\le2(C+K)$.

Put $\iota_{\rm root}=1$ if a root face is present, and $0$ otherwise.
\begin{lemma}[Vertex budget]\label{lem:budget}
$\sum_jE_j=\Dnorm$ and $\sum_jV_j\le\Dnorm+\iota_{\rm root}$; hence $\sum_j(V_j-E_j)\le\iota_{\rm root}$.
\end{lemma}
\begin{proof}
The first identity counts every edge according to its multiplicity. For the second, $\sum_jV_j=\sum_z\max_{e\ni z}m(e)$. Couple the incidence of each arriving facet at $z$ to that of the next departing facet in its trace, using cyclic order for non-root traces. All incidences are paired except the first and last incidences of the root trace. Let $k(z)$ count these exceptions at $z$, so $\sum_zk(z)=2\iota_{\rm root}$. Coupled incidences lie on distinct edges because the traces do not backtrack. Hence
\[
2\max_{e\ni z}m(e)\le\sum_{e\ni z}m(e)+k(z).
\]
Summing over vertices gives $2\sum_jV_j\le2\Dnorm+2\iota_{\rm root}$.
\end{proof}

\subsection{Labelled embeddings}

A \emph{map of labelled graphs} is a graph map that commutes with $e\mapsto\inv e$ and preserves labels. Since $p_a=p_{\inv a}$, the frequency $p_{\ell(e)}$ of the label of an edge does not depend on its orientation.

\begin{lemma}[Labelled embeddings]\label{lem:embed}
Let $C'>0$, and let $H$ be a graph without loops, multiple edges or isolated vertices, with $V$ vertices and $E\le C'L$ edges, labelled by $\ell_H$, with a distinguished root if $\iota=1$ and no root if $\iota=0$, where $\iota\in\{0,1\}$. Let $Y\subseteq V(H)$ contain the root when present, meet every component and contain every vertex of degree other than $2$; suppressing the vertices outside $Y$ gives a multigraph with $Q$ edges, which we orient. For $v\notin Y$ let $t_v,u_v$ be the labels of the edges entering and leaving $v$ along its path, in the chosen orientation. Then for large $n$ the $\mu_n$-expected number of maps of labelled graphs $H\to\Gamma_\Phi$ that are injective on vertices and map the root, when present, to $x_0$ is at most
\[
e^{o(L)}\,\pmin^{-Q}\,n^{V-E-\iota}\prod_{v\notin Y}w(t_v,u_v),
\]
where $o(L)$ depends only on $\pi$, $c$ and $C'$.
\end{lemma}

\begin{proof}
As $\Gamma_\Phi$ is immersed, a map is determined by its injective vertex map $g$, and exists if and only if the $E$ edges of $H$, transported by $g$, are edges of $\Gamma_\Phi$. These are $E$ distinct prescriptions, as $g$ is injective and $H$ has no multiple edges, so by Lemma~\ref{lem:cond-presc}(b) they hold with probability at most $e^{o(L)}\prod_\epsilon(np_{\ell_H(\epsilon)})^{-1}$, the product being over the edges $\epsilon$ of $H$. For $v\notin Y$, the vertex $g(v)$ must have outlabels $\inv{t_v}$ and $u_v$. If $u_v=\inv{t_v}$, no injective map exists, since the two edges at $v$ would have the same image; otherwise $g(v)$ ranges over $V_{\inv{t_v}}\cap V_{u_v}$, a set of $np_{\inv{t_v},u_v}=np_{t_v}w(t_v,u_v)$ vertices. The root has one possible image, while the other vertices of $Y$ have at most $n$ images each. Pair each $v\notin Y$ with the edge entering it along its path, whose factor $(np_{t_v})^{-1}$ cancels $np_{t_v}$. The $Q$ unpaired edges, the last edges of the paths, contribute at most $(n\pmin)^{-Q}$. As $V-|Y|=E-Q$, the bound follows.
\end{proof}

\subsection{Pieces}

Orient each chain; we regard the orientations as part of the pattern, which multiplies the bound of Lemma~\ref{lem:count} by at most $2^{c_1/2}$. Identify a chain $\chain$ with $[0,|\chain|]$, its edges sitting in order at the half-integers $\frac12,\frac32,\dots,|\chain|-\frac12$. By Lemma~\ref{lem:thin}(c), each facet $\xi$ lies over the edge at a position $x(\xi)$ of a chain $\chain(\xi)$, within a traversal of direction $\theta(\xi)=\pm1$ relative to the orientation of $\chain(\xi)$. For a labelling $\ell$ of $\Delta$ let $\ell_{\chain}(x)$ be the label of the edge at position $x$ of $\chain$, oriented along $\chain$; with $a^{+1}:=a$ and $a^{-1}:=\inv a$, we have $\ell(\lambda(\xi))=\ell_{\chain(\xi)}(x(\xi))^{\theta(\xi)}$.

A $T$-labelling $\ell$ of $\Delta$ is \emph{admissible} if it is reduced (distinct oriented edges with the same initial vertex carry distinct letters) and $\ell(\lambda(\xi'))=\inv{\ell(\lambda(\xi))}$ for all paired facets $\xi,\xi'$. If $\psi$ realizes the pattern in $\Gamma_\Phi$, the pullback of the edge labelling of $\Gamma_\Phi$ is admissible: it is reduced because $\psi$ is injective and $\Gamma_\Phi$ is immersed, and paired facets receive inverse letters by the definition of a realization.

Cut each interior segment at the vertices where, on either side, the trace passes from one traversal to the next. The resulting \emph{elementary segments} have, on each side, an interval of facets within a single traversal, and these two intervals are paired in reverse order. Each of the at most $c_2$ passages between consecutive traversals lies in the interior of at most one interval, so by Lemma~\ref{lem:segments}(c) there are at most $g:=\lceil I+K+c_2\rceil$ elementary segments.

\begin{lemma}[Comparison maps]\label{lem:comparison}
Let $\alpha$ and $\alpha'$ be the two intervals of an elementary segment, over chains $\chain$ and $\chain'$. There are $\sigma\in\{\pm1\}$ and $a\in\Z$ such that $x(\xi')=\sigma x(\xi)+a$ whenever $\xi\in\alpha$ is paired with $\xi'\in\alpha'$, and every admissible labelling satisfies $\ell_{\chain'}(\sigma x+a)=\ell_{\chain}(x)^\sigma$ for $x=x(\xi)$, $\xi\in\alpha$. If $\chain'=\chain$ and $\sigma=1$, then $a\ne0$.
\end{lemma}

\begin{proof}
Let $\theta,\theta'$ be the directions of the traversals containing $\alpha,\alpha'$. Moving to the next facet of $\alpha$ changes $x$ by $\theta$, while the partner moves to the previous facet of $\alpha'$ and its position changes by $-\theta'$. Hence $x(\xi')=\sigma x(\xi)+a$ with $\sigma:=-\theta\theta'$ and $a\in\Z$. Admissibility gives $\ell_{\chain'}(x')^{\theta'}=\ell_{\chain}(x)^{-\theta}$, that is, $\ell_{\chain'}(x')=\ell_{\chain}(x)^\sigma$. If $\chain'=\chain$, $\sigma=1$ and $a=0$, paired facets lie over the same edge, contrary to Lemma~\ref{lem:pattern-props}(b).
\end{proof}

The inverse of $x\mapsto\sigma x+a$ is $x\mapsto\sigma x-\sigma a$.

\medskip\noindent\textbf{Piece length.} Let $a_1,\dots,a_g$ be the integers $a$ of the elementary segments, padded with zeros, and put $N:=\lfloor L^{1/(2g+2)}\rfloor$ and $R:=\lfloor N^{1/2}\rfloor$. By the pigeonhole principle applied to the points $(\{ka_i/L\})_{i\le g}\in[0,1)^g$, $0\le k\le N^g$, where $\{\cdot\}$ denotes the fractional part, some integer $1\le u\le N^g$ satisfies $\dist(ua_i/L,\Z)\le1/N$ for all $i$; take the least one and put $s:=L/(uR)$. Since $\dist(Ry,\Z)\le R\dist(y,\Z)$ for integers $R$,
\begin{equation}\label{eq:s}
LN^{-g-1/2}\le s\le L/R,\qquad\dist(a_i/s,\Z)\le R/N\quad(1\le i\le g).
\end{equation}
Thus $s\ge L^{1/2}$, $s/L\le1/R$ and $R/N\le N^{-1/2}$, so $\kappa:=R/N+2/s\to0$ and $s/L\to0$ as $n\to\infty$, uniformly over patterns. We take $n$ so large that $s\ge8$, $(1-\kappa)s>(s+1)/2$, and the estimates~\eqref{eq:cost-bounds}--\eqref{eq:cost-count} below apply.

\medskip\noindent\textbf{Pieces.} For a chain $\chain$ and an integer $k\ge0$ with $(k+1)s\le|\chain|$, the \emph{piece} $(\chain,k)$ is the subpath of $\chain$ formed by the edges at positions in $[ks,(k+1)s)$; the rest of $\chain$, its \emph{end fragment}, has fewer than $s+1$ edges. A piece $\beta$ has $h_\beta\in\{\lfloor s\rfloor,\lceil s\rceil\}$ edges, and we identify it with its set of positions; its multiplicity is $m(\beta):=m(\chain)$. An \emph{occurrence} of $\beta$ is the set of facets over $\beta$ within one traversal of $\chain$, an interval of one star. There are $m(\beta)$ occurrences of $\beta$, distinct occurrences are disjoint, and since each has at least $s-1$ facets there are at most $2CL/s$ of them in all. For a labelling $\ell$, the \emph{word} of $\beta$ is $v_\beta:=\ell_{\chain}(x_1)\cdots\ell_{\chain}(x_{h_\beta})$, where $x_1<\dots<x_{h_\beta}$ are the positions of $\beta$; we write $v_\beta(x)$ for its letter at position $x$.

\begin{lemma}[Alignment]\label{lem:align}
Let $x\mapsto\sigma x+a$ be the map of an elementary segment or its inverse, and let $\beta=(\chain,k)$ be a piece. Choose $\nu\in\Z$ with $|a-\nu s|\le(R/N)s$, and put $k':=k+\nu$ if $\sigma=1$ and $k':=\nu-k-1$ if $\sigma=-1$. Then at least $(1-\kappa)s$ consecutive positions $x$ of $\beta$ satisfy $\sigma x+a\in[k's,(k'+1)s)$.
\end{lemma}

\begin{proof}
$\nu$ exists by~\eqref{eq:s}, since the shift of the inverse map is $\pm a$. Write $a=\nu s+\varrho$. Apart from the choice of open and closed endpoints when $\sigma=-1$, the image of $[ks,(k+1)s)$ is the translate of $[k's,(k'+1)s)$ by $\varrho$. The two intervals overlap in length $s-|\varrho|\ge(1-R/N)s$. Allowing for their endpoints, this overlap contains at least $(1-R/N)s-2=(1-\kappa)s$ half-integers. These are images of consecutive positions of $\beta$, since $a\in\Z$.
\end{proof}

\subsection{Joins and the constraint graph}

Two distinct occurrences $\omega,\omega'$ are \emph{joined} if some elementary segment pairs at least $(1-\kappa)s$ consecutive facets of $\omega$ with facets of $\omega'$.

\begin{lemma}[Joins]\label{lem:joins}
\begin{enumerate}[(a)]
\item Each occurrence is joined to at most one other.
\item The number $b_\star$ of occurrences not joined to any other occurrence satisfies $s\,b_\star\le \operatorname{free}(D)+o(L)$.
\end{enumerate}
\end{lemma}

\begin{proof}
(a) An occurrence has fewer than $s+1<2(1-\kappa)s$ facets, so two sets of at least $(1-\kappa)s$ of its facets share a facet; its unique partner lies in both partner occurrences, which therefore coincide.

(b) Call the first and last facets of the intervals of elementary segments \emph{endpoints}; there are at most $4g$. Call an occurrence in a star $f$ \emph{exceptional} if an endpoint in $f$ is within $3s$ positions, in the linear order of $f$, of one of its facets. Occurrences are disjoint intervals of at least $s-1$ facets, so each endpoint makes at most $9$ occurrences exceptional, and at most $36g$ are exceptional.

Let $\omega$ be a non-exceptional occurrence of $\beta=(\chain,k)$ containing a paired facet, which lies in an interval $\alpha$ of an elementary segment with map $x\mapsto\sigma x+a$ to the other interval $\alpha'$, over $\chain'$. Then $\omega\subseteq\alpha$, and all facets of $\omega$ are more than $3s$ facets away from the ends of $\alpha$. By Lemma~\ref{lem:align}, at least $(1-\kappa)s$ consecutive facets of $\omega$ have partners over positions in $[k's,(k'+1)s)$. These partners are more than $3s$ facets away from the ends of $\alpha'$, which lies in one traversal of $\chain'$, so their positions lie in $(3s,|\chain'|-3s)$. Hence $k'\ge0$ and $(k'+1)s\le|\chain'|$, so $(\chain',k')$ is a piece, the partners lie in its occurrence $\omega'$ in that traversal, and $\omega'\ne\omega$ as $\alpha\cap\alpha'=\emptyset$. So $\omega$ is joined to $\omega'$.

Consequently, a non-exceptional unjoined occurrence consists of free facets, and there are at most $\operatorname{free}(D)/(s-1)$ such occurrences. As $\operatorname{free}(D)\le CL$ and $s/(s-1)\le1+2/s$, we get $s\,b_\star\le36gs+\operatorname{free}(D)+2CL/s=\operatorname{free}(D)+o(L)$.
\end{proof}

\noindent\textbf{The constraint graph.} Let $\cC$ be the multigraph whose vertices are the pieces, with an edge, or \emph{link}, joining the pieces of $\omega$ and $\omega'$ for each joined pair $\{\omega,\omega'\}$; a link from a piece to itself is a \emph{self-link}. For each link choose an elementary segment realizing the join, let $\beta$ be the piece of $\omega$, and let $J\subseteq\beta$ be the set of positions under the at least $(1-\kappa)s$ consecutive facets of $\omega$ that it pairs with facets of $\omega'$. The link records $J$ and the restriction to $J$ of the corresponding comparison map, which takes values in the piece $\beta'$ of $\omega'$.

A \emph{piece assignment} on a set $U$ of pieces is a family $(v_\beta)_{\beta\in U}$ of reduced words with $|v_\beta|=h_\beta$. It is \emph{compatible} if $v_{\beta'}(\sigma x+a)=v_\beta(x)^\sigma$ for every link inside $U$, with recorded data $(J,x\mapsto\sigma x+a)$, and every $x\in J$. The piece words of an admissible labelling are reduced, being labels of non-backtracking paths, and by Lemma~\ref{lem:comparison} they are compatible on every set of pieces.

\begin{lemma}[Propagation]\label{lem:propagation}
Let $U$ be a component of $\cC$. In this lemma, the error term $o(s)$ denotes $(\kappa s+1)\log|T|$, which is sublinear in $s$ because $\kappa\to0$ and $s\to\infty$.
\begin{enumerate}[(a)]
\item For $\beta_0\in U$, each reduced word of length $h_{\beta_0}$ extends in at most $e^{o(s)|U|}$ ways to a compatible assignment on $U$.
\item If $U$ contains a self-link, there are at most $e^{o(s)|U|}$ compatible assignments on $U$.
\end{enumerate}
\end{lemma}

\begin{proof}
(a) Along a spanning tree of $U$ rooted at $\beta_0$, the word of a parent determines at least $(1-\kappa)s$ letters of the child word through the link, leaving fewer than $\kappa s+1$ free letters.

(b) Let the self-link at $\beta$ have data $(J,x\mapsto\sigma x+a)$, and put $h:=h_\beta$. If $\sigma=1$, then $a\ne0$ and $v_\beta(x+a)=v_\beta(x)$ for $x\in J$; the graph on the positions of $\beta$ with edges $\{x,x+a\}$, $x\in J$, is a forest with $h-|J|<\kappa s+1$ components on which $v_\beta$ is constant, so there are at most $|T|^{\kappa s+1}$ choices for $v_\beta$. If $\sigma=-1$, then $v_\beta(a-x)=\inv{v_\beta(x)}$ for $x\in J$. The $x\in J$ with $a-x\in J$ form a nonempty interval symmetric about $a/2$, since $2|J|-h\ge2(1-\kappa)s-(s+1)>0$. It contains $a/2$ if $a$ is odd, giving a letter equal to its inverse, and $a/2\pm\frac12$ if $a$ is even, giving two adjacent mutually inverse letters; so no compatible assignment exists. In the translation case, apply (a) with $\beta_0=\beta$; in the reflection case, there are no compatible assignments.
\end{proof}

\begin{lemma}[Multiplicities]\label{lem:mult}
For a component $U$ of $\cC$, let $S_U$ and $z_U$ be the sum and the maximum of $m(\beta)$ over $\beta\in U$, and let $b_U$ be the number of unjoined occurrences of pieces of $U$. If $U$ has no self-link, then $2z_U\le S_U+b_U$.
\end{lemma}

\begin{proof}
Let $m(\beta)=z_U$. At least $z_U-b_U$ occurrences of $\beta$ are joined, by Lemma~\ref{lem:joins}(a) to pairwise distinct occurrences, which belong to pieces of $U$ other than $\beta$. So $z_U-b_U\le S_U-z_U$.
\end{proof}

\subsection{Counting realizations by stages}

The embedding estimate charges a piece word $v$ a factor $P(v)$, while
propagation lets us count the words on a whole component of the constraint
graph after choosing one of them. We compare these two contributions by
recording each positive word weight to within a factor of $e$. For $P(v)>0$ put
$d(v):=\lfloor-\log P(v)\rfloor$. If $|v|=h\in(s-1,s+1)$, then, uniformly
for all sufficiently large $n$,
\begin{equation}\label{eq:cost-bounds}
 \frac{\delta s}{2}\le d(v)\le a_0s,
 \qquad e^{-d(v)-1}<P(v)\le e^{-d(v)}.
\end{equation}
Indeed, $P(v)^2\le\sum_{|u|=h}P(u)^2\le e^{-2\delta h}$ gives
$d(v)\ge\delta h-1\ge\delta s/2$ for large $s$, whereas the $h-1<s$
positive factors of $P(v)$ are at least $\wmin=e^{-a_0}$.
These bounds are uniform because $s\ge L^{1/2}\to\infty$.

Fix an oriented pattern from Lemma~\ref{lem:count}, and fix one admissible
choice of its mesh $s$ from the construction above. Fix also an integer $d_\beta$
in $[\delta s/2,a_0s]$ for each piece $\beta$. Throughout the following
estimates, a labelling is required to be admissible and to satisfy
$P(v_\beta)>0$ and $d(v_\beta)=d_\beta$ for every piece. Let
$\mathcal E_{\mathbf d}$ be the event that the pattern has a realization with
such a labelling, where $\mathbf d=(d_\beta)_\beta$.

For $1\le j\le m_\star$, let $X_j$ count the graph maps
$\Delta_j\to\Gamma_\Phi$ that are injective on vertices and whose pulled-back
labelling extends to a labelling of the whole pattern satisfying these
requirements. When $j=1$ and a root is present, require the root to map to
$x_0$. A realization counted by $\mathcal E_{\mathbf d}$ restricts to a map
counted by each $X_j$. Thus Markov's inequality gives
\begin{equation}\label{eq:markov}
 \mu_n(\mathcal E_{\mathbf d})\le\min_j\mathbb E X_j.
\end{equation}

\begin{lemma}[Stage bound]\label{lem:stage}
For the pattern and integers fixed above, $\mathbb E X_j\le B_j$, where
\[
\log B_j:=(V_j-E_j-\mathbf1_{j=1}\iota_{\rm root})\log n
-\sum_{m(\beta)\ge j}d_\beta
+\sum_{U\in\cU_j}(2d_U-\delta s)+o(L).
\]
Here $\cU_j$ consists of the components $U$ of $\cC$ having no self-link and
$z_U\ge j$, and $d_U:=\min_{\beta\in U}d_\beta$. The error is uniform over
the patterns and the choices of integers.
\end{lemma}

\begin{proof}
First fix the pulled-back labelling $\ell_j$ of $\Delta_j$. Apply
Lemma~\ref{lem:embed} with $H=\Delta_j$, $Y=V(\Delta_j)\cap V_\star$,
$Q\le c_1/2$, $C'=C$ and $\iota=\mathbf1_{j=1}\iota_{\rm root}$; the root,
when required, belongs to $V_\star$. Each chain in $\Delta_j$ inherits its orientation from $\Delta$, so the turn factors inside a piece $\beta$ multiply to $P(v_\beta)$. Discard the factors outside the pieces of multiplicity at least $j$. By~\eqref{eq:cost-bounds}, the remaining product is at most $\exp(-\sum_{m(\beta)\ge j}d_\beta)$. The expected number
of maps with this labelling is therefore at most
\[
 e^{o(L)}n^{V_j-E_j-\mathbf1_{j=1}\iota_{\rm root}}
 \exp\Bigl(-\sum_{m(\beta)\ge j}d_\beta\Bigr).
\]

It remains to count the possible labellings $\ell_j$. The end fragments have
at most $|T|^{(s+1)c_1/2}=e^{o(L)}$ labellings. In a component $U$ of $\cC$
with $z_U<j$, no piece word appears at this stage. If $z_U\ge j$ and $U$ has
a self-link, Lemma~\ref{lem:propagation}(b) gives at most $e^{o(s)|U|}$
compatible assignments. Otherwise choose a piece $\beta_0\in U$ with
$d_{\beta_0}=d_U$. Its length is fixed. For any such length $h$ and integer $d$, each
word with $d(v)=d$ contributes more than $e^{-2d-2}$ to the squared-weight sum. Hence
\begin{equation}\label{eq:cost-count}
 \#\{v: |v|=h,\ P(v)>0,\ d(v)=d\}
 \le e^{2d+2}\sum_{|v|=h}P(v)^2
 \le e^{2d+2-2\delta h}\le e^{2d-\delta s}.
\end{equation}
The final inequality uses $h>s-1$ and holds uniformly for large $s$.
Thus there are at most $e^{2d_U-\delta s}$ choices for the word on $\beta_0$.
Lemma~\ref{lem:propagation}(a) gives at most $e^{o(s)|U|}$ extensions.
This is a valid upper bound even if $m(\beta_0)<j$: by the definition of
$X_j$, each stage labelling extends to a compatible assignment on all of $U$,
and counting those assignments can only overcount their restrictions.
Since $\sum_U|U|\le2CL/s$, all propagation errors together contribute
$e^{o(L)}$, uniformly. Multiplication gives the stated bound.
\end{proof}

\begin{lemma}[Summation]\label{lem:sum}
For every fixed pattern and choice $\mathbf d$ as above,
\[
 \sum_{j=1}^{m_\star}\log B_j
 \le-\frac\delta2\Dnorm+a_0\,\operatorname{free}(D)+o(L),
\]
with a uniform error.
\end{lemma}

\begin{proof}
Lemma~\ref{lem:budget} gives
$\sum_j(V_j-E_j)-\iota_{\rm root}\le0$, so the powers of $n$ make no
positive contribution. There are at most $2(C+K)$ stages, and their uniform
errors sum to $o(L)$. A piece $\beta$ occurs in the stages
$j\le m(\beta)$, whereas $U\in\cU_j$ for $j\le z_U$. Hence the contribution
of a component $U$ is
\[
 F_U=-\sum_{\beta\in U}m(\beta)d_\beta
       +\begin{cases}
          0,&U\text{ has a self-link},\\
          z_U(2d_U-\delta s),&U\text{ has no self-link}.
        \end{cases}
\]
We claim $F_U\le-\frac\delta2sS_U+a_0sb_U$. For a self-link component,
this follows directly from $d_\beta\ge\delta s/2$. Otherwise
\[
 F_U\le(2z_U-S_U)d_U-\delta sz_U.
\]
If $2z_U\le S_U$, the lower bound $d_U\ge\delta s/2$ gives
$F_U\le-\delta sS_U/2$. If $2z_U>S_U$, use $d_U\le a_0s$,
$2z_U-S_U\le b_U$ from Lemma~\ref{lem:mult}, and $z_U>S_U/2$ to obtain
the claimed bound.

Finally, $\sum_US_U$ counts occurrences. Each has at most $s+1$ facets,
and together they cover all but at most $c_2(s+1)$ facets, so
\[
 s\sum_US_U\ge\Dnorm-c_2(s+1)-\sum_US_U=\Dnorm-o(L).
\]
Lemma~\ref{lem:joins}(b) gives
$s\sum_Ub_U=s\,b_\star\le\operatorname{free}(D)+o(L)$. Summing the component
bounds proves the result.
\end{proof}

\begin{proof}[Proof of Theorem~\ref{thm:local}]
A violating diagram determines an oriented pattern and an admissible
labelling of its trace graph. Every piece word has positive weight: consecutive
letters $a,b$ of a non-backtracking trace satisfy $b\ne\inv a$ and
$p_{\inv a,b}>0$. Thus~\eqref{eq:cost-bounds} assigns an allowed integer
$d_\beta$ to every piece, and the configuration lies in
$\mathcal E_{\mathbf d}$. The conditions
$L\le\Dnorm\le CL$ and $\operatorname{free}(D)\le\varepsilon\Dnorm$
depend only on the pattern class, by Lemma~\ref{lem:pattern-props}(c).

There are at most $2^{c_1/2}L^{c_3}$ oriented pattern classes. For each one
we have fixed its mesh $s$, and, for all sufficiently large $n$,
$L^{1/2}\le s\le L$. It has at most $2CL/s\le2C\sqrt L$ pieces, and each
integer $d_\beta$ has at most $a_0L+2$ choices. Consequently the total number
of pairs consisting of an oriented pattern and a choice $\mathbf d$ is at most
\begin{equation}\label{eq:cost-profile-count}
 2^{c_1/2}L^{c_3}(a_0L+2)^{2C\sqrt L}
 =\exp\bigl(O(\sqrt L\log L)\bigr)=e^{o(L)}.
\end{equation}
This bound is uniform even though $s$ may vary between patterns.

Fix one such choice satisfying the two violation inequalities. Since
$a_0\varepsilon\le\delta/16$, Lemma~\ref{lem:sum} gives, uniformly,
\[
 \sum_j\log B_j\le-\frac{7\delta}{16}\Dnorm+o(L)
 \le-\frac\delta4\Dnorm
\]
for large $n$, using $\Dnorm\ge L$. There are at most
$M:=\lceil2(C+K)\rceil$ stages, so one of them has
$\log B_j\le-\delta\Dnorm/(4M)\le-\delta L/(4M)$. By~\eqref{eq:markov}
and Lemma~\ref{lem:stage},
\[
 \mu_n(\mathcal E_{\mathbf d})\le e^{-\delta L/(4M)}.
\]
Summing this bound over~\eqref{eq:cost-profile-count}, the probability of
any violating diagram is at most $e^{o(L)-\delta L/(4M)}$, which tends to zero.
\end{proof}

\begin{remark}
Besides (H1) and (H2), the proof uses only the pair counts $|V_{\inv a}\cap V_b|=np_aw(a,b)$, the positive lower bound $\wmin$ for nonzero turn weights, and the decay $\sum_{|v|=h}P(v)^2\le e^{-2\delta h}$.
\end{remark}

\section{A linear isoperimetric inequality}\label{sec:global}

Mghirbi and Douglas obtain the area bound for the characteristic-two examples by separating cone pictures into small planar pieces \cite[Lemma~3.7 and Proposition~4.4]{md2026}. The local estimate of Theorem~\ref{thm:local} concerns triangular diagrams with only boundedly many stars. Here we pass from that estimate to all disk diagrams.

We use Gromov's local-to-global principle, in the combinatorial form of Papasoglu~\cite{Pap96} and Ollivier~\cite[App.~A]{ollivier} (Lemma~\ref{lem:localglobal}). It requires a class of disks closed under taking subdisks, with a uniform bound on face lengths. Stars may be arbitrarily long, so we first cut each long star along chords into a chain of tiles with bounded perimeter. The first and last edges of each chord are chosen so that the boundary of any union of consecutive tiles traces a closed, cyclically non-backtracking path in $\Gamma$ (Section~\ref{ss:chains}). We then reassemble each tiled subdisk with boundedly many tiles into a reduced diagram over $\cP_\Phi$, to which Theorem~\ref{thm:local} applies (Section~\ref{ss:tiled}).

Let $\varepsilon$ be as in~\eqref{eq:eps} and $\Dcl$ as in Lemma~\ref{lem:closure}, and put
\begin{equation}\label{eq:epsone}
\lambda:=2\Dcl+8,\qquad\kappa:=\Bigl\lceil\max\Bigl\{10^{10},\ \frac{3\cdot10^7\lambda}{\varepsilon}\Bigr\}\Bigr\rceil,\qquad\varepsilon_1:=\frac{1}{2\cdot10^4\,\kappa(\lambda+1)}.
\end{equation}
These constants depend only on $\pi$ and $c$.

\begin{theorem}[Linear isoperimetric inequality]\label{thm:global}
Assume \textup{(H1)--(H3)} and let $\varepsilon_1$ be as in~\eqref{eq:epsone}. With high probability, every reduced diagram $D$ over $\cP_\Phi$ that is either $0$-hole (possibly filamentous) or spherical, has boundary words over $T$, and has at least one face satisfies $|\bd D|>\varepsilon_1\Dnorm$. In particular, w.h.p.\ no spherical diagram over $\cP_\Phi$ with a face is reduced.
\end{theorem}

The proof is given in Section~\ref{ss:globalproof}.

\subsection{The local-to-global principle}\label{ss:localglobal}

Let $X$ be a simply connected combinatorial $2$-complex. A \emph{disk drawn in $X$} is a combinatorial map $\phi\colon E\to X$, where $E$ is a finite combinatorial $2$-complex homeomorphic to a closed disk and $\phi$ maps each open cell of $E$ homeomorphically onto an open cell of $X$. Its \emph{length} $L_c(\phi)$ is the number of edges of $\bd E$, and its \emph{area} is $A_c(\phi):=\sum_t|\bd t|^2$, summed over the faces $t$ of $E$, where $|\bd t|$ is the length of the boundary path of $t$. A \emph{subdisk} of $E$ is a subcomplex homeomorphic to a closed disk; it is a union of closed faces of $E$, and $\phi$ restricts to a disk drawn in $X$.

\begin{lemma}[{Papasoglu, Ollivier~\cite[App.~A, Prop.~42]{ollivier}}]\label{lem:localglobal}
Let $X$ be a simply connected combinatorial $2$-complex whose faces have at most $\ell$ edges, and let $\mathsf P$ be a property of disks drawn in $X$ such that the restriction of a disk with $\mathsf P$ to any subdisk has $\mathsf P$. Suppose that for some integer $K\ge10^{10}\ell$, every disk $\phi$ with $\mathsf P$ and $K^2/4\le A_c(\phi)\le480K^2$ satisfies $L_c(\phi)^2\ge2\cdot10^{14}A_c(\phi)$. Then every disk $\phi$ with $\mathsf P$ and $A_c(\phi)\ge K^2$ satisfies $L_c(\phi)\ge A_c(\phi)/(10^4K)$.
\end{lemma}

This is Ollivier's extension to polygonal complexes of Papasoglu's theorem~\cite{Pap96}. As noted in~\cite[App.~A]{ollivier}, the proof only passes to subdisks of a given disk drawn in $X$, which is why it applies to any property $\mathsf P$ as above.

\subsection{Arcs and disk components}\label{ss:arcs}

\begin{lemma}[Arcs]\label{lem:arcs}
Let $D$ be a $0$-hole, possibly filamentous, or spherical $T$-diagram with $N\ge1$ stars. Then $D$ has fewer than $3N$ arcs.
\end{lemma}

\begin{proof}
Let $D^\circ$ be the union of the closed faces of $D$. Every face vertex other than an apex lies in $D_T$: it is an endpoint of the $T$-edge of a triangle, or the initial point of the root star boundary. The $Z$-edges at an apex lie in its star. Thus Lemma~\ref{lem:stars} gives $D^\circ$ a coarser CW structure with $1$-skeleton $D_T$ and the $N$ stars as $2$-cells. Consequently $\chi(D^\circ)=V-E+N$, where $V$ and $E$ count the vertices and edges of $D_T$.

If $D$ is spherical, then $D^\circ=S^2$. Otherwise view $D$ in $S^2=\R^2\cup\{\infty\}$. The complement $S^2\setminus D^\circ$ is connected and nonempty: $S^2\setminus D$ is connected, and each point of $D\setminus D^\circ$ lies on no face and has a connected neighbourhood disjoint from $D^\circ$ that meets $S^2\setminus D$. Alexander duality~\cite[Cor.~3.45]{Hat02} gives $\tilde H^1(D^\circ)=H^2(D^\circ)=0$, so $\chi(D^\circ)\ge1$. In both cases, taking degrees in $D_T$, we obtain
\begin{equation}\label{eq:euler}
\sum_v(\deg v-2)=2E-2V=2N-2\chi(D^\circ)\le2N-2.
\end{equation}

Every vertex $v$ of $D_T$ of degree one must be the unmarked vertex of the root face. Indeed, let $a$ be its $T$-edge and $c$ a triangle on it. The adjacent corner across the $Z$-edge at $v$ is again a face corner. If that face is a triangle, its other edge must be $a$, giving consecutive triangles of a star whose traces have inverse letters. Immersion forces a backtrack, contrary to Lemma~\ref{lem:stars}. Thus the adjacent face is the root face. There is at most one such vertex. Writing $u\le1$ for the number of degree-one vertices, equation~\eqref{eq:euler} gives
\[
\sum_{\deg v\ne2}\deg v\le3(2N-2+u)+u\le6N-2<6N.
\]
Each arc has two ends counted in this sum, so there are fewer than $3N$ arcs.
\end{proof}

\begin{lemma}[Disk components]\label{lem:diskcomp}
Let $Y$ be a $0$-hole, possibly filamentous, diagram over some presentation, with at least one face. Call two faces of $Y$ \emph{adjacent} if they share an edge, and let $C_1,\dots,C_r$ be the unions of the closed faces in the classes of the equivalence relation generated by adjacency. Then each $C_j$ is a subcomplex of $Y$ homeomorphic to a closed disk, each face of $Y$ lies in exactly one $C_j$, and $\sum_j|\bd C_j|\le|\bd Y|$, where $|\bd C_j|$ is the number of edges of the boundary circle of $C_j$.
\end{lemma}

\begin{proof}
We regard the planar embedding $Y\subset\R^2$ as piecewise linear. Put $S^2=\R^2\cup\{\infty\}$ and $\Omega:=S^2\setminus Y$, which is connected because $Y$ has no holes. Fix $j$, let $\mathcal F_j$ be the class whose union is $C_j$, and note that $C_j$ is connected. The edge-ends at a vertex $v$ cut a small disk about $v$ into sectors, each of which lies in an open face or in $\Omega$. If two consecutive sectors at $v$ lie in faces $c$ and $c'$, both faces contain the edge between them, so $c=c'$ or $c$ and $c'$ are adjacent. Hence each maximal run of consecutive face sectors consists of corners from a single adjacency class.

We claim that the corners at $v$ of faces of $\mathcal F_j$ form one such run, or all sectors at $v$. Otherwise there are corners $s,s'$ at $v$ of faces $c,c'\in\mathcal F_j$ and sectors $\sigma,\sigma'\subseteq\Omega$ at $v$ such that $s$ and $s'$ separate $\sigma$ from $\sigma'$ in the cyclic order around $v$. Choose distinct faces $c=c_0,c_1,\dots,c_p=c'$ in $\mathcal F_j$, consecutive ones sharing an edge $\eta_i$, and a piecewise linear simple closed curve $\beta$ that leaves $v$ into $s$, runs through the open faces $c_0,\dots,c_p$ in turn, crossing each $\eta_i$ once at an interior point, and returns to $v$ through $s'$. Near $v$ the two germs of $\beta$ separate $\sigma$ from $\sigma'$, so $\sigma$ and $\sigma'$ lie in different components of $S^2\setminus\beta$; but they lie in the connected set $\Omega$, which is disjoint from $\beta$. This proves the claim.

Consequently every vertex of $C_j$ has a neighbourhood in $C_j$ homeomorphic to a disk or a half-disk. So does every interior point of an edge $\eta$ of $C_j$: a face of $\mathcal F_j$ lies on one side of $\eta$, and on the other side lies either $\Omega$ or a face containing $\eta$, which is adjacent to the first or equal to it, hence in $\mathcal F_j$. Points of open faces are interior. Thus $C_j$ is a compact connected surface, possibly with boundary. Each boundary edge has a face of $\mathcal F_j$ on one side and $\Omega$ on the other.

The complement $S^2\setminus C_j$ is connected. Otherwise, as $\Omega$ is connected, some component $W$ of $S^2\setminus C_j$ is disjoint from $\Omega$. The frontier of $W$ lies in $C_j$ but not in an open face, so it lies in the $1$-skeleton of $C_j$; and it is infinite, for otherwise $W$ would be open and closed in the connected set $S^2\setminus\operatorname{Fr}W$, hence equal to it, and would meet $C_j$. So the frontier contains an interior point $p$ of an edge $\eta$ of $C_j$. Near $p$, the side of $\eta$ that meets $W$ lies in a face of $\mathcal F_j$ or in $\Omega$, by the previous paragraph; the first is impossible as $W\cap C_j=\emptyset$, and the second as $W\cap\Omega=\emptyset$.

A compact connected surface without boundary embedded in $S^2$ must be all of $S^2$, by invariance of domain. Since $C_j$ is disjoint from the nonempty set $\Omega$, it has a boundary circle. By the Schoenflies theorem (e.g.~\cite{Moi77}), each boundary circle $J$ of $C_j$ is a piecewise linear Jordan curve, and $C_j\setminus J$, being connected, lies in one component of $S^2\setminus J$; the other component is then a component of $S^2\setminus C_j$ with frontier $J$. As $S^2\setminus C_j$ is connected, $C_j$ has exactly one boundary circle $J$, and $C_j$ is the complement of the open disk bounded by $J$ on the other side, that is, a closed disk.

Finally, the classes partition the faces of $Y$. An edge of $\bd C_j$ has $\Omega$ on exactly one side, so it is read once in the boundary word of $Y$, and it lies on no face outside $\mathcal F_j$, so it is not a boundary edge of another $C_{j'}$. Hence $\sum_j|\bd C_j|\le|\bd Y|$.
\end{proof}

\subsection{Chords and chain refinements}\label{ss:chains}

In this subsection $\Gamma$ is as in Lemma~\ref{lem:closure}. For a vertex $z$ let $\mathrm{Out}(z)$ be the set of oriented edges with initial vertex $z$; since $\Gamma$ has no loops, $|\mathrm{Out}(z)|=\deg z\ge3$.

\begin{lemma}[Chords]\label{lem:chords}
Let $u$ and $v$ lie in the same connected component of $\Gamma$, and let $A\subseteq\mathrm{Out}(u)$ and $B\subseteq\mathrm{Out}(v)$ have at most two elements each. There is a non-backtracking path $\rho$ from $u$ to $v$ with $1\le|\rho|\le\Dcl L+2$ whose first edge is not in $A$ and whose last edge is not the reverse of an element of $B$.
\end{lemma}

\begin{proof}
Put $k:=\lceil\log_2n\rceil+2$. For a vertex $z$ and $h\in\mathrm{Out}(z)$, the proof of Lemma~\ref{lem:closure} gives a closed non-backtracking path $J_z(h)$ at $z$ of length at most $2k$ whose first edge is $h$ and whose last edge is $\inv h$. Choose $f\in\mathrm{Out}(u)\setminus A$ and $g\in\mathrm{Out}(v)\setminus B$, and let $u'$ and $v'$ be the terminal vertices of $f$ and $g$. If $u'\ne v'$, let $R$ be a geodesic from $u'$ to $v'$, with first edge $r_1$ and last edge $r_2$; choose $h_1\in\mathrm{Out}(u')\setminus\{\inv f,r_1\}$ and $h_2\in\mathrm{Out}(v')\setminus\{\inv g,\inv{r_2}\}$, and put
\[
\rho:=f\,J_{u'}(h_1)\,R\,J_{v'}(h_2)\,\inv g.
\]
If $u'=v'$, choose $h\in\mathrm{Out}(u')\setminus\{\inv f,\inv g\}$ and put $\rho:=f\,J_{u'}(h)\,\inv g$. Each constituent path is non-backtracking. At each junction the outgoing edge is not the reverse of the incoming one: $h_1\ne\inv f$, $r_1\ne h_1$, $h_2\ne\inv{r_2}$ and $\inv g\ne h_2$ (respectively $h\ne\inv f$ and $\inv g\ne h$). So $\rho$ is non-backtracking; its first edge is $f\notin A$, and its last edge is $\inv g$, the reverse of $g\notin B$. Finally, $L\ge\frac c2\log n$ gives $4k\le\frac8{c\log2}L+12$, so $|\rho|\le2+4k+d_0L\le\Dcl L+2$.
\end{proof}

Let $D$ be a reduced $0$-hole, possibly filamentous, or spherical diagram over $\cP_\Phi$ with boundary words over $T$, and let $f=\St(o)$ be a star of $D$ with $h:=|\bd f|$ facets, a chosen first face, facets $\xi_1,\dots,\xi_h$ and lift $\gamma_f=e_1\cdots e_h$, where $e_i:=\gamma_f(\xi_i)$. Indices are taken modulo $h$, and $x_i$ denotes the terminal vertex of $e_i$, which is the initial vertex of $e_{i+1}$. We identify $Q_o$ (Lemma~\ref{lem:stars}(c)) with a regular $h$-gon in $\R^2$ whose vertices $v_0,\dots,v_{h-1}$ occur in positive order, $\xi_i$ being the edge from $v_{i-1}$ to $v_i$; position $i$ on $\bd Q_o$ means the vertex $v_i$. Put
\[
F(i):=\{\inv{e_i},\,e_{i+1}\}\subseteq\mathrm{Out}(x_i).
\]
A \emph{chord from position $i$ to position $i'$} is a non-backtracking path from $x_i$ to $x_{i'}$ whose first edge is not in $F(i)$ and whose last edge is not the reverse of an element of $F(i')$.

\begin{definition}[Chain refinement of a star]\label{def:chain}
If $h\le\lambda L$, the \emph{tiles} of $f$ consist of the single tile $T_0:=Q_o$, and $r:=0$. If $h>\lambda L$, put $r:=\lceil h/(4L)\rceil$, $q_i:=\lfloor(i-1)h/(2r)\rfloor$ for $1\le i\le2r$, and $a_k:=q_k$, $b_k:=q_{2r+1-k}$ for $1\le k\le r$. For each $k$, the vertices $x_{a_k}$ and $x_{b_k}$ lie on the same trace $\gamma_f$. Lemma~\ref{lem:chords}, with $A=F(a_k)$ and $B=F(b_k)$, therefore gives a chord $\rho_k$ from position $a_k$ to position $b_k$ with $|\rho_k|\le\Dcl L+2$. Let $d_k$ be the straight segment from $v_{a_k}$ to $v_{b_k}$ in $Q_o$, subdivided into $|\rho_k|$ edges labelled so that $d_k$, read from $v_{a_k}$, spells $\ell(\rho_k)$. The \emph{tiles} of $f$ are: $T_0$, bounded by the arc of $\bd Q_o$ from $v_{b_1}$ to $v_{a_1}$ and by $d_1$; $T_k$ for $1\le k\le r-1$, bounded by $d_k$, $d_{k+1}$ and the arcs from $v_{a_k}$ to $v_{a_{k+1}}$ and from $v_{b_{k+1}}$ to $v_{b_k}$; and $T_r$, bounded by $d_r$ and the arc from $v_{a_r}$ to $v_{b_r}$. Here an arc from $v_i$ to $v_{i'}$ consists of the edges $\xi_{i+1},\dots,\xi_{i'}$ of $\bd Q_o$.
\end{definition}
For an interval $I=\{k_1,\dots,k_2\}$ in $\{0,\dots,r\}$, let $Q_I$ be the union of the tiles $T_k$, $k\in I$, and let $\pi_I$ be the closed path in $\Gamma$ obtained by reading $\bd Q_I$ positively, lifting each edge $\xi_i$ to $e_i$ and each chord $d_k$, read from $v_{a_k}$, to $\rho_k$.

\begin{lemma}[Chains]\label{lem:chains}
Let $f$ be a star refined as in Definition~\ref{def:chain}, and let $L\ge3$.
\begin{enumerate}[(a)]
\item If $h>\lambda L$, then $0=a_1<\dots<a_r<b_r<\dots<b_1<h$, each arc of $\bd Q_o$ between consecutive points of $\{v_{a_k},v_{b_k}\}$ has between $L$ and $2L$ edges, the segments $d_k$ are pairwise disjoint, and the tiles are convex polygons such that $T_{k-1}\cap T_k=d_k$ for $1\le k\le r$ and $T_k\cap T_{k'}=\emptyset$ for $|k-k'|\ge2$.
\item Each facet of $f$ lies on exactly one tile, each tile carries at least $L$ facets of $f$, and each tile has at most $\lambda L$ edges.
\item For every interval $I$, $Q_I$ is a convex polygon, and $\pi_I$ is a closed, cyclically non-backtracking path that traverses every facet of $f$ lying on a tile of $I$ and has length at most $\sum_{k\in I}|\bd T_k|$.
\item In particular, the boundary of each tile $T_k$, read positively, spells the label of the closed, cyclically non-backtracking path $\pi_{\{k\}}$, of length at most $\lambda L$.
\end{enumerate}
\end{lemma}

\begin{proof}
(a) Since $h>\lambda L\ge8L$, we have $r\ge2$, and $r\ge h/(4L)$ and $r<h/(4L)+1$ give
\[
\tfrac43L\le\frac{2Lh}{h+4L}<\frac h{2r}\le2L.
\]
The differences $q_{i+1}-q_i$ for $1\le i<2r$, and $h-q_{2r}+q_1$, lie in $\{\lfloor h/(2r)\rfloor,\lceil h/(2r)\rceil\}$; as $h/(2r)>\frac43L\ge L+1$ for $L\ge3$, and $2L$ is an integer, they lie between $L$ and $2L$. In particular the $q_i$ increase from $q_1=0$ to $q_{2r}<h$, which gives the order of the $a_k$ and $b_k$, as $a_r=q_r<q_{r+1}=b_r$. Two segments joining distinct vertices of a convex polygon are disjoint if their four endpoints do not interleave in the cyclic order. Thus the $d_k$ are pairwise disjoint and cut $Q_o$ into the $r+1$ convex polygons $T_0,\dots,T_r$, with the stated intersections.

(b) The arcs partition the edges of $\bd Q_o$, which are the facets of $f$, and each arc lies on exactly one tile. By (a), $T_0$ and $T_r$ carry one arc each and the other tiles two; so each tile carries at least $L$ facets of $f$ and has at most $4L+2(\Dcl L+2)\le\lambda L$ edges. If $h\le\lambda L$, the single tile carries all $h\ge L$ facets and has $h\le\lambda L$ edges.

(c) If $I=\{0,\dots,r\}$, then $Q_I=Q_o$ and $\pi_I=\gamma_f$. Otherwise $Q_I$ is the part of $Q_o$ between the chords $d_{k_1}$ (if $k_1\ge1$) and $d_{k_2+1}$ (if $k_2\le r-1$), a convex polygon, and $\bd Q_I$, read positively, alternates between arcs of $\bd Q_o$, read forwards, and these chords: $d_{k_1}$ is read from $v_{b_{k_1}}$ to $v_{a_{k_1}}$ and $d_{k_2+1}$ from $v_{a_{k_2+1}}$ to $v_{b_{k_2+1}}$. The arcs lift to subpaths of the cyclically non-backtracking path $\gamma_f$, and the chords lift to $\rho_k$ or its reverse. These paths are all non-backtracking. It remains to check the junctions of $\pi_I$, each of which lies at a chord endpoint $v_i$ between an arc and a chord. The arc ending at $v_i$ has last edge $e_i$, and the arc starting there has first edge $e_{i+1}$. By the definition of a chord, its edge directed away from $x_i$ lies outside $F(i)=\{\inv{e_i},e_{i+1}\}$, regardless of which endpoint of the chord is at $v_i$. This excludes a backtrack both when an arc enters the chord and when the chord enters an arc. Thus $\pi_I$ is cyclically non-backtracking. It traverses the arcs of all tiles of $I$, and each edge of $\bd Q_I$ is an edge of exactly one tile of $I$.

(d) Take $I=\{k\}$ in (c) and use (b).
\end{proof}

\subsection{Tiled disks}\label{ss:tiled}

Let $\cR_\lambda$ be the set of labels of the closed, cyclically non-backtracking paths of length at most $\lambda L$ in $\Gamma$, a finite set closed under cyclic permutation and inversion.

Let $\cD$ be the set of reduced diagrams over $\cP_\Phi$ that are either $0$-hole (possibly filamentous) or spherical, have boundary words over $T$, and have at least one face. A \emph{chain refinement} $D^\natural$ of $D\in\cD$ is obtained by choosing tiles for every star $f=\St(o)$ as in Definition~\ref{def:chain} and replacing the cells of $\St(o)$ by the $q_o$-images of the open tiles and of the open edges and interior vertices of the segments $d_k$. By Lemma~\ref{lem:stars}(c) this is a cell structure on the same space as $D$; its $2$-cells are the tiles, the tile $T_k$ of $f$ being attached by $q_o|_{\bd T_k}$. Label the chord edges as in Definition~\ref{def:chain}, mark each tile at a vertex, and orient it positively. By Lemma~\ref{lem:chains}(d), this makes $D^\natural$ a diagram over $\langle S\mid\cR_\lambda\rangle$. If $D$ is $0$-hole, then so is $D^\natural$, with the same boundary word and possibly with filaments; if $D$ is spherical, so is $D^\natural$. Every facet of $D$ lies on exactly one tile (Lemma~\ref{lem:chains}(b)), and a chord edge lies on exactly the two consecutive tiles of one star that contain it. For a union $E$ of tiles, put $A(E):=\sum_{t\subseteq E}|\bd t|$. In particular, $\Dnorm\le A(D^\natural)$.

A \emph{tiled disk} is a subcomplex $E$ of a chain refinement $D^\natural$ of some $D\in\cD$ that is homeomorphic to a closed disk. It is a union of closed tiles, since a subcomplex homeomorphic to a closed disk is the union of the closed $2$-cells it contains; if $D$ is spherical, then $E\ne S^2$, and we view $E$ in $\R^2$ by stereographic projection from a point of $S^2\setminus E$. Thus $E$ is a diagram over $\langle S\mid\cR_\lambda\rangle$, and a subdisk of a tiled disk is a tiled disk.

\begin{lemma}[Reassembly]\label{lem:reassembly}
Let $M\ge1$ be an integer. On the event of Theorem~\ref{thm:local} for $K=M$, $I=3M$ and $C=\lambda M$, every tiled disk $E$ with at most $M$ tiles satisfies $|\bd E|>(\varepsilon/\lambda)A(E)$.
\end{lemma}

\begin{proof}
Let $E\subseteq D^\natural$ with $D\in\cD$. For a star $f=\St(o)$ of $D$, the set of indices $k$ with $T_k\subseteq E$ is a disjoint union of maximal intervals $I$ of $\{0,\dots,r\}$; the corresponding sets of tiles are the \emph{pieces} of $E$, and a piece $P$ has the convex polygon $Q_P:=Q_I\subseteq Q_o$ and the closed path $\pi_P:=\pi_I$. Two consecutive tiles of $f$ lying in $E$ belong to the same piece. Every tile of $E$ lies in exactly one piece, so $E$ has at most $M$ pieces.

The open polygon $\operatorname{int}Q_P$ consists of the open tiles of $P$ and of the open edges and interior vertices of the chords between consecutive tiles of $P$, all of which lie in $E$. As $q_o$ maps $\operatorname{int}Q_o$ homeomorphically onto the open set $\St(o)$, it maps $\operatorname{int}Q_P$ homeomorphically onto an open subset $U_P$ of $E$, and the sets $U_P$ for distinct pieces are disjoint. In each $Q_P$ choose an interior point $o_P$ and join it by straight segments to all vertices of the cell structure of $\bd Q_P$, that is, the vertices $v_i$ in $Q_P$ and the subdivision vertices of the chords on $\bd Q_P$; transporting this cell structure by $q_o$ to $U_P$, for every piece $P$, gives a cell structure $E^*$ on $E$, whose $2$-cells are the triangles spanned by the points $o_P$ and the edges of the circles $\bd Q_P$. Label the segment from $o_P$ to a vertex of $\bd Q_P$ at which $\pi_P$ visits the vertex $x$ of $\Gamma$ by $z_x$ towards $o_P$, and let the triangle over an edge of $\bd Q_P$ that traces to the edge $e$ of $\pi_P$ bear $r_{e^+}$, marked at $o_P$ and oriented so that it lies over $e$. Then $E^*$ is a $0$-hole, non-filamentous diagram over $\cP_\Phi$ with the same boundary word as $E$, a word over $T$.

We next check that $E^*$ is reduced. Two triangles sharing a $Z$-edge are consecutive triangles of one fan and lie over consecutive edges $e,e'$ of $\pi_P$. If they formed a reducible pair, they would bear the same relator. Since $r_{e^+}$ determines $e^+$, this would give $e'\in\{e,\inv e\}$. The equality $e'=e$ is impossible because $\Gamma$ has no loops, and $e'=\inv e$ is excluded by the non-backtracking property of $\pi_P$ (Lemma~\ref{lem:chains}(c)). Let $\eta$ be a $T$-edge of $E^*$ lying on two triangles. It is not a chord edge: a chord edge lies only on two consecutive tiles of one star, and if both lie in $E$ they belong to one piece and the chord edge lies in some $U_P$. So $\eta$ is a $T$-edge of $D$ lying on two faces $c,c'$ of $D$, and the triangles of $E^*$ on its two sides lie over $e_c$ and $e_{c'}$, which are distinct even up to orientation by Lemma~\ref{lem:stars}(d); so they bear different relators.

By Lemma~\ref{lem:stars}, the stars of $E^*$ are the fans about the points $o_P$, with traces $\pi_P$; so $E^*$ has at most $M$ stars, fewer than $3M$ arcs by Lemma~\ref{lem:arcs}, and $\norm{E^*}=\sum_P|\pi_P|$. If a piece $P$ contains $t_P$ tiles, then $\pi_P$ traverses at least $L$ facets of $f$ from each tile. Lemma~\ref{lem:chains}(b),(c) therefore gives
\[
t_PL\le|\pi_P|\le\sum_{t\in P}|\bd t|\le t_P\lambda L\le\lambda|\pi_P|.
\] Hence $L\le\norm{E^*}\le A(E)\le\lambda ML$ and $A(E)\le\lambda\norm{E^*}$. Theorem~\ref{thm:local} gives $|\bd E|=|\bd E^*|>\varepsilon\norm{E^*}\ge(\varepsilon/\lambda)A(E)$.
\end{proof}

\begin{proposition}\label{prop:tiled}
Put $\ell:=\lceil\lambda L\rceil$, $K_0:=\kappa\ell$ and $M:=\lceil480\kappa^2(\lambda+1)^2\rceil$, with $\lambda,\kappa$ as in~\eqref{eq:epsone}. On the event of Theorem~\ref{thm:local} for $K=M$, $I=3M$ and $C=\lambda M$, every tiled disk $E$ satisfies $|\bd E|\ge2\varepsilon_1A(E)$.
\end{proposition}

\begin{proof}
Let $X_\lambda$ be the Cayley complex of $\langle S\mid\cR_\lambda\rangle$, the universal cover of its presentation complex. It is a simply connected combinatorial $2$-complex whose faces have at most $\lambda L\le\ell$ edges. A tiled disk $E$ is simply connected and each tile boundary spells a relator, so mapping a vertex of $E$ to a vertex of $X_\lambda$ extends to a label-preserving combinatorial map $\phi_E\colon E\to X_\lambda$~\cite[Ch.~V]{LS77}, a disk drawn in $X_\lambda$ with $L_c(\phi_E)=|\bd E|$ and $A_c(\phi_E)=\sum_t|\bd t|^2$. Let $\mathsf P$ be the property of being of the form $\phi_E$ for a tiled disk $E$. It passes to subdisks, since a subdisk of a tiled disk is a tiled disk and the restriction of a label-preserving map is label-preserving.

Every tile of $E$ has between $L$ and $\lambda L$ edges, by Lemma~\ref{lem:closed-nb} and Lemma~\ref{lem:chains}(d); so $L\,A(E)\le A_c(\phi_E)\le\lambda L\,A(E)$, and $E$ has at most $A_c(\phi_E)/L^2$ tiles. If $K_0^2/4\le A_c(\phi_E)\le480K_0^2$, then $E$ has at most $480\kappa^2\ell^2/L^2\le M$ tiles, as $\ell\le(\lambda+1)L$; so Lemma~\ref{lem:reassembly} and $K_0\ge\kappa\lambda L$ give
\begin{align*}
|\bd E|^2&>\Bigl(\frac\varepsilon\lambda A(E)\Bigr)^2\ge\Bigl(\frac\varepsilon{\lambda^2}\Bigr)^2\frac{A_c(\phi_E)}{L^2}\,A_c(\phi_E)\\
&\ge\Bigl(\frac\varepsilon{\lambda^2}\Bigr)^2\frac{K_0^2}{4L^2}\,A_c(\phi_E)\ge\Bigl(\frac{\varepsilon\kappa}{2\lambda}\Bigr)^2A_c(\phi_E)\ge2\cdot10^{14}A_c(\phi_E),
\end{align*}
since $\varepsilon\kappa/(2\lambda)\ge1.5\cdot10^7$. As $K_0\ge10^{10}\ell$, Lemma~\ref{lem:localglobal} applies with $K=K_0$: if $A_c(\phi_E)\ge K_0^2$, then
\[
|\bd E|\ge\frac{A_c(\phi_E)}{10^4K_0}\ge\frac{L\,A(E)}{10^4\kappa(\lambda+1)L}=2\varepsilon_1A(E).
\]
If $A_c(\phi_E)<K_0^2$, then $E$ has fewer than $K_0^2/L^2\le M$ tiles, and Lemma~\ref{lem:reassembly} gives $|\bd E|>(\varepsilon/\lambda)A(E)\ge2\varepsilon_1A(E)$, as $10^4\kappa(\lambda+1)\ge\lambda/\varepsilon$.
\end{proof}

\subsection{Proof of Theorem~\ref{thm:global}}\label{ss:globalproof}

\begin{proof}
Let $M$ be as in Proposition~\ref{prop:tiled} and $N_0:=\lceil\lambda(1+1/(2\varepsilon_1))\rceil$. Let $E_n$ be the intersection of the events of Lemma~\ref{lem:expansion} and of Theorem~\ref{thm:local} for $(K,I,C)=(M,3M,\lambda M)$ and for $(K,I,C)=(N_0,3N_0,N_0)$; then $\mu_n(E_n)\to1$. Let $n$ be large, so that $L\ge\max\{3,\frac c2\log n\}$, and let $\Phi\in E_n$. By Lemmas~\ref{lem:immersion} and~\ref{lem:expansion}, (H2) and $\girth(\Gamma_\Phi)\ge L\ge3$, the graph $\Gamma_\Phi$ satisfies the hypotheses of Lemma~\ref{lem:closure}, so every $D\in\cD$ has chain refinements and Proposition~\ref{prop:tiled} applies. Let $D\in\cD$ and fix a chain refinement $D^\natural$.

\emph{The $0$-hole case.} Then $D^\natural$ is a $0$-hole, possibly filamentous, diagram with the same boundary word as $D$ and at least one face. Let $C_1,\dots,C_r$ be as in Lemma~\ref{lem:diskcomp} for $Y=D^\natural$. They are tiled disks, and every tile lies in exactly one of them; so Lemma~\ref{lem:diskcomp} and Proposition~\ref{prop:tiled} give
\[
|\bd D|\ge\sum_j|\bd C_j|\ge2\varepsilon_1\sum_jA(C_j)=2\varepsilon_1A(D^\natural)\ge2\varepsilon_1\Dnorm>\varepsilon_1\Dnorm.
\]

\emph{The spherical case.} We show that this is impossible. Let $t_0$ be a tile of $D^\natural$ and $Y:=S^2\setminus\operatorname{int}t_0$, where $\operatorname{int}t_0$ is the open $2$-cell, an open disk. By Alexander duality~\cite[Cor.~3.45]{Hat02}, the subcomplex $Y$ is connected. Projecting from a point of $\operatorname{int}t_0$, $Y$ with the cells and labels of $D^\natural$ other than $t_0$ is a $0$-hole, possibly filamentous, diagram over $\langle S\mid\cR_\lambda\rangle$, whose boundary word is read along $\bd t_0$: an edge with $t_0$ on one side is read once, and an edge with $t_0$ on both sides is a filament of $Y$, read twice. So $|\bd Y|=|\bd t_0|\le\lambda L$. If $Y$ has no face, then $A(D^\natural)=|\bd t_0|$. Otherwise let $C_1,\dots,C_r$ be as in Lemma~\ref{lem:diskcomp} for $Y$; they are tiled disks, and every tile other than $t_0$ lies in exactly one of them, so
\[
\lambda L\ge|\bd Y|\ge\sum_j|\bd C_j|\ge2\varepsilon_1\sum_jA(C_j)=2\varepsilon_1\bigl(A(D^\natural)-|\bd t_0|\bigr).
\]
In both cases $\Dnorm\le A(D^\natural)\le\lambda L(1+1/(2\varepsilon_1))\le N_0L$. Every star has at least $L$ facets, so $D$ has at most $N_0$ stars, and it has fewer than $3N_0$ arcs by Lemma~\ref{lem:arcs}. As $L\le\Dnorm\le N_0L$, Theorem~\ref{thm:local} gives $\operatorname{free}(D)>\varepsilon\Dnorm>0$, whereas $\operatorname{free}(D)=0$ because $D$ is spherical.
\end{proof}

\section{Proof of the outlabel criterion}\label{sec:proofs}

We now turn the two isoperimetric inequalities into the four assertions of Theorem~\ref{intro:construction}. The global inequality gives hyperbolicity and asphericity. Root separation needs one further application of the local inequality, to a sphere containing one root face.

\subsection{The cone complex and the auxiliary presentation}

The vertices of $X_\Phi$ are the vertex $*$ of $R_S$ and a cone point $o_{\Gamma'}$ for each component $\Gamma'$ of $\Gamma_\Phi$. Its edges are the petals of $R_S$ and, for every vertex $x$ of $\Gamma_\Phi$, a cone edge $\zeta_x$ from $*$ to the cone point of the component of $x$. An edge $e$ from $x$ to $x'$ with label $a$ contributes a triangle attached along $a\zeta_{x'}\inv{\zeta_x}$.

\begin{proposition}\label{prop:cone}
Let $\Phi\in\Omega_n$, define $\cP_\Phi$ by the formulas of Section~\ref{sec:diagrams}, and let $X_{\cP}$ be its presentation complex and $G_{\cP}$ its fundamental group.
\begin{enumerate}[(a)]
\item For a word $v$ over $T$, the following statements are equivalent: $v=1$ in $G_\Phi$; $v=1$ in $G_{\cP}$; and $v$ is null-homotopic in $X_\Phi$, viewed as a loop in $R_S$. Moreover,
\[
\pi_1(X_\Phi,*)\cong G_\Phi=\langle S\mid B_\Phi\rangle.
\]
\item The natural maps $\pi_j(X_\Phi)\to\pi_j(X_{\cP})$ are injective for all $j\ge1$.
\item If a word $v$ over $T$ bounds a $0$-hole, possibly filamentous, diagram over $\cP_\Phi$ with $N$ faces, then in $F(S)$ it is a product of at most $N$ conjugates of elements of $B_\Phi\cup B_\Phi^{-1}$.
\end{enumerate}
\end{proposition}

\begin{proof}
Choose a spanning tree and a base vertex in each component of $\Gamma_\Phi$, and let $p_x$ be the tree path from its base vertex to $x$. Define $\Upsilon\colon F(S\cup Z)\to F(S)$ to fix $S$ and send $z_x$ to $\ell(p_x)^{-1}$. For an edge $e$ from $x$ to $x'$,
\[
\Upsilon(r_e)=\ell(p_xe\inv{p_{x'}}).
\]
This is the label of a fundamental cycle, or its inverse, if $e$ is outside the spanning trees, and is freely trivial otherwise. For a closed path $\gamma=e_1\cdots e_h$ based at $u$, cancellation of successive $Z$-letters gives
\begin{equation}\label{eq:telescope}
z_u r_{e_1}\cdots r_{e_h}z_u^{-1}=\ell(\gamma)
\quad\text{in }F(S\cup Z).
\end{equation}
Part (c) follows by applying $\Upsilon$ to the expression supplied by the van Kampen lemma.

For the topological comparison, add to $X_\Phi$ one extra edge from $*$ to $o_{\Gamma'}$ for each component $\Gamma'$, obtaining $X^+$. These extra edges form a tree. Collapsing it is a homotopy equivalence~\cite[Prop.~0.17]{Hat02}; the quotient has one vertex, one loop for every letter of $S\cup Z$, and the triangles of $\cP_\Phi$. Thus the quotient is $X_{\cP}$, with the petal loops unchanged. Sending each added edge to any cone edge in the same component defines a retraction $X^+\to X_\Phi$. Consequently $X_\Phi\to X^+$, and hence $X_\Phi\to X_{\cP}$, is injective on every homotopy group. This proves (b).

If $v=1$ in $G_\Phi$, it is a product of conjugates of labels of closed paths and their inverses; by~\eqref{eq:telescope} it is trivial in $G_{\cP}$. If $v$ is trivial in $G_{\cP}$, the homotopy equivalence $X^+\to X_{\cP}$ and the retraction show that it is null-homotopic in $X_\Phi$. Conversely, a null-homotopy in $X_\Phi$ gives one in $X_{\cP}$. Part (c) then puts $v$ in the normal closure of $B_\Phi$, so $v=1$ in $G_\Phi$.

Applying $\Upsilon$ to~\eqref{eq:telescope} also shows that the labels of all closed paths lie in the normal closure of $B_\Phi$. Hence $G_\Phi=\langle S\mid B_\Phi\rangle$. Finally, the petals generate $\pi_1(X_\Phi,*)$: an edge loop is a product of petals and of paths $\zeta_x\inv{\zeta_y}$ through a cone point, where $x$ and $y$ lie in the same component. The triangles identify such a path with the label of a path from $x$ to $y$ in that component. The equivalence already proved identifies the kernel, giving (a).
\end{proof}

\begin{lemma}[Gersten]\label{lem:gersten}
Assume $\girth(\Gamma_\Phi)\ge3$. If every spherical diagram over $\cP_\Phi$ with at least one face has a reducible pair, then $\pi_2(X_{\cP})=0$.
\end{lemma}

\begin{proof}
The three letters of each relator of $\cP_\Phi$ are distinct. Thus a reducible pair is precisely a pair of faces mapped to the same $2$-cell as mirror images across their common edge. The hypothesis says that $X_{\cP}$ is diagrammatically reducible. Diagrammatically reducible complexes are aspherical~\cite{Ger87}.
\end{proof}

\subsection{Hyperbolicity and the absence of short relations}

\begin{proof}[Proof of parts (1) and (2) of Theorem~\ref{intro:construction}]
Take $n$ large and $\Phi$ in the event of Theorem~\ref{thm:global}. If $v=1$ in $G_\Phi$, Proposition~\ref{prop:cone}(a) and the van Kampen lemma give a reduced $0$-hole diagram $D$ over $\cP_\Phi$ with boundary word $v$. Every face has one facet. Thus, when $D$ has a face, the global inequality gives
\[
|D|=\Dnorm<|v|/\varepsilon_1.
\]
Proposition~\ref{prop:cone}(c) gives the required expression as a product of conjugates of relators. A finite presentation with a linear isoperimetric inequality defines a hyperbolic group~\cite[Ch.~III.H, Thm.~2.9]{bh}.

If $v$ is nonempty and cyclically reduced, $D$ cannot be a tree, since the boundary word of a tree is freely trivial. Thus $D$ has a face and hence a star. Since $D$ has no root face, that star has at least $L$ facets by Lemma~\ref{lem:stars}(b). Hence
\[
|v|=|\bd D|>\varepsilon_1\Dnorm\ge\varepsilon_1L.
\]
In particular, once $\varepsilon_1L\ge1$, no generator is trivial.
\end{proof}

\subsection{Asphericity, torsion-freeness and non-elementarity}

\begin{proof}[Proof of part (3) of Theorem~\ref{intro:construction}]
Theorem~\ref{thm:global} excludes reduced spherical diagrams over $\cP_\Phi$ with a face. Lemma~\ref{lem:gersten} gives $\pi_2(X_{\cP})=0$, and Proposition~\ref{prop:cone}(b) gives $\pi_2(X_\Phi)=0$. The universal cover $\widetilde X_\Phi$ is a simply connected $2$-complex with $H_2(\widetilde X_\Phi)=0$. The Hurewicz and Whitehead theorems imply that it is contractible. Thus $X_\Phi$ is a finite two-dimensional $K(G_\Phi,1)$.

Its cellular chain complex gives a free $\mathbb ZG_\Phi$-resolution of $\mathbb Z$ of length two. A subgroup of prime order $q$ would inherit such a resolution by restriction. This would force $H^3(\mathbb Z/q;\mathbb F_q)=0$, contrary to the cohomology of a cyclic group of prime order~\cite[Ch.~VIII, \S2]{brown}. Hence $G_\Phi$ is torsion-free.

Parts (1) and (2) show that $G_\Phi$ is hyperbolic and nontrivial. Since it is torsion-free, it is infinite. It is not virtually cyclic. Indeed, $X_\Phi$ has $1+\#\pi_0(\Gamma_\Phi)$ vertices, $r+n$ edges and $|E(\Gamma_\Phi)|$ triangles. Since $2|E(\Gamma_\Phi)|\ge3n$,
\[
\chi(X_\Phi)=1+\#\pi_0(\Gamma_\Phi)-r-n+|E(\Gamma_\Phi)|
\ge2-r+n/2>0
\]
for large $n$. An infinite cyclic subgroup of finite index $j$ would give a finite $j$-sheeted cover that is a $K(\mathbb Z,1)$, hence homotopy equivalent to a circle~\cite[Thm.~1B.8]{Hat02}. Its Euler characteristic would be both zero and $j\chi(X_\Phi)>0$, a contradiction.
\end{proof}

\begin{remark}
The isoperimetric constant $1/\varepsilon_1$ depends only on $\pi$ and $c$. The fundamental-cycle relators have length at least $L$ and, when breadth-first spanning trees are used, length $O(L)$ by Lemma~\ref{lem:expansion}(b). Accordingly, this proof does not give a hyperbolicity constant independent of $n$.
\end{remark}

\subsection{Root separation}

\begin{proof}[Proof of part (4) of Theorem~\ref{intro:construction}]
Let $d_0$ be the component-diameter constant of Lemma~\ref{lem:expansion}, and put
\[
K:=\lfloor d_0(1+1/\varepsilon_1)\rfloor+1.
\]
Take $n$ large and intersect the events of Lemma~\ref{lem:expansion} and Theorem~\ref{thm:global} with the event of Theorem~\ref{thm:local} for the prescribed root $x_0$ and parameters $(K,I,C)=(K,3K,K)$. Their intersection has probability tending to one. Suppose that the label of some path from $x_0$ to $y\ne x_0$ is trivial in $G_\Phi$.

Choose a geodesic $\gamma=e_1\cdots e_d$ from $x_0$ to $y$. Then $1\le d\le d_0L$, and $\ell(\gamma)=1$ in $G_\Phi$, because the labels of any two paths with these endpoints differ by a closed-path relation. The global inequality supplies a reduced diagram for $\ell(\gamma)$ over $\cP_\Phi$ with at most $d/\varepsilon_1$ faces. Moreover,
\[
z_{x_0}^{-1}\ell(\gamma)z_y=r_{e_1}\cdots r_{e_d}.
\]
Consequently
\[
z_y^{-1}z_{x_0}
=\bigl(z_y^{-1}\ell(\gamma)z_y\bigr)
  \bigl(r_{e_1}\cdots r_{e_d}\bigr)^{-1}
\]
is a product of at most $d(1+1/\varepsilon_1)$ conjugates of relators of $\cP_\Phi$ and their inverses. Choose a reduced $0$-hole diagram $D$ with this boundary word and at most $d_0(1+1/\varepsilon_1)L$ faces.

We first check that $D$ is a closed disk, so that it can be capped by a root face. A filament would be traversed twice in the boundary word of $D$. Since that word has length two, it would then be freely trivial, contrary to $y\ne x_0$. Thus there are no filaments, and the boundary consists of two distinct edges. Apply Lemma~\ref{lem:diskcomp} to $D$. Each boundary edge of a disk component lies on a face on only one side, so it is one of these two edges, and the sum of the component perimeters is at most two. No component can have perimeter one: its boundary would read $z_x^{\pm1}$, whereas the homomorphism $F(S\cup Z)\to\mathbb Z$ sending every letter of $S$ to zero and every $z_x$ to one kills all the relators of $\cP_\Phi$ but does not kill $z_x$. There is therefore exactly one disk component. Every edge lies on a face, so this component is all of $D$.

By the Schoenflies theorem, the closure of the complement of $D$ in $S^2$ is a disk. Make it one face bearing the root relator $z_y^{-1}z_{x_0}$, with its marking and orientation chosen accordingly. The resulting spherical diagram $\widehat D$ is reduced: $D$ was reduced, and the new face cannot form a reducible pair because no other face bears a root relator. Thus $\widehat D$ is a $T$-diagram with exactly one root face, and
\[
\norm{\widehat D}\le d_0(1+1/\varepsilon_1)L<KL.
\]

There must be a star other than the root star $f_0=\St(o)$. By Lemma~\ref{lem:stars}(b), the trace of $f_0$ is a non-backtracking path of positive length. Since $\Gamma_\Phi$ is immersed, its label, read along the closed path $\bd_o$ in $\widehat D$, is a nonempty reduced word. If $f_0$ were the only star, the complement $S^2\setminus f_0$ would be the graph $\widehat D_T$. The star is an open disk by Lemma~\ref{lem:stars}(c). Alexander duality then implies that $\widehat D_T$ is connected and has vanishing first cohomology, hence is a tree. But every closed path in a tree reads a freely trivial word, a contradiction.

Every star other than $f_0$ has at least $L$ facets. Therefore
\[
L\le\norm{\widehat D}<KL,
\]
there are at most $K$ stars, and Lemma~\ref{lem:arcs} gives fewer than $3K$ arcs. The local inequality now yields
\[
\operatorname{free}(\widehat D)>\varepsilon\norm{\widehat D}>0.
\]
This is impossible for a spherical diagram. Root separation follows.
\end{proof}

\section{Several independent parts}\label{sec:partitioned}

The two-sided zero-divisor construction uses independent matchings on two prescribed parts of the vertex set. We retain this partition when sampling: using only the averaged link distribution would allow edges between the parts and define a different random graph model.

\subsection{The partitioned model}

\begin{definition}[Partitioned model]\label{def:partitioned}
Fix an integer $J\ge1$, balanced link distributions $\pi_1,\ldots,\pi_J$, each using every letter, and positive rational weights $\nu_1,\ldots,\nu_J$ with sum one. Write $\boldsymbol\pi=(\pi_i,\nu_i)_{i\le J}$. A positive integer $n$ is admissible when every $n\nu_i\pi_i(\tau)$ is an integer. Partition $V_n$ into sets $V^{(i)}$ of size $\nu_i n$, and prescribe links so that
\[
|\{x\in V^{(i)}:\tau_n(x)=\tau\}|=n\nu_i\pi_i(\tau).
\]
For $V_a^{(i)}:=V_a\cap V^{(i)}$, define
\[
\Omega_n^{\mathrm{pt}}
:=\{\Phi\in\Omega_n:\phi_s(V_s^{(i)})=V_{s^{-1}}^{(i)}
                  \text{ for every }i,s\}
\cong\prod_{i=1}^J\prod_{s\in S}\Bij(V_s^{(i)},V_{s^{-1}}^{(i)}).
\]
The overall proportions of the links are given by $\overline\pi=\sum_i\nu_i\pi_i$. The graph, group and cone complex are defined as before, and $\Gamma_\Phi$ is the disjoint union of its subgraphs $\Gamma_\Phi^{(i)}$ on the parts.
\end{definition}

The unconditioned measure is uniform on $\Omega_n^{\mathrm{pt}}$. Put
\[
\begin{gathered}
\pmin:=\min_{i,a}\nu_i p_a(\pi_i),\qquad
\dmax:=\max_i\dmax(\pi_i),\\
c_*(\boldsymbol\pi):=\min\left\{\frac1{\log(2|T|/\pmin)},\frac1{2\log\dmax}\right\}.
\end{gathered}
\]
For $0<c<c_*(\boldsymbol\pi)$ and $L=\lfloor c\log n\rfloor$, let $\mu_n$ be the uniform measure on configurations in $\Omega_n^{\mathrm{pt}}$ with girth at least $L$. The girth condition holds precisely when it holds in every part. The conditioned law is therefore still a product, so the parts remain independent.

\begin{theorem}[Outlabel criterion for several parts]\label{thm:partitioned}
Assume that every $\pi_i$ satisfies \textup{(H1)--(H3)}, let $0<c<c_*(\boldsymbol\pi)$, and fix a root $x_0\in V_n$ for each admissible $n$. There is $\varepsilon_1=\varepsilon_1(\boldsymbol\pi,c)>0$ such that, with $\mu_n$-probability tending to one as $n$ tends to infinity through admissible values:
\begin{enumerate}[(1)]
\item every word $v$ over $T$ trivial in $G_\Phi$ is a product in $F(S)$ of at most $|v|/\varepsilon_1$ conjugates of elements of $B_\Phi\cup B_\Phi^{-1}$;
\item no nonempty cyclically reduced word of length at most $\varepsilon_1L$ is trivial in $G_\Phi$;
\item $X_\Phi$ is a finite two-dimensional $K(G_\Phi,1)$, and $G_\Phi$ is torsion-free and non-elementary hyperbolic;
\item the label of every path from $x_0$ to a different vertex is nontrivial in $G_\Phi$.
\end{enumerate}
The same conclusions hold, with the same $\varepsilon_1$, if the girth threshold is replaced by integers $L_n$ satisfying $c\log n-C_0\le L_n\le c\log n$ for a fixed $C_0\ge0$; assertion \textup{(2)} then uses $L_n$.
\end{theorem}

The theorem allows any fixed finite number of parts, with different link distributions and different positive proportions. The collision bound is checked separately in each part. The two-part zero-divisor construction is one example.

\begin{example}[The two-sided zero-divisor construction]\label{ex:partitioned-zero}
Take the distributions $\pi_A,\pi_B$ of Example~\ref{ex:openai}, two parts of equal size, and $n=2n'$ with $n'$ divisible by $2(q^2+q+1)^2$. These are exactly the vertex sets and matchings of~\cite[\S2]{seed}. The girth threshold in that construction is $L_n=\lfloor c\log n'\rfloor$, where $c=1/100$, and it satisfies $c\log n-2\le L_n\le c\log n$. Here $\pmin=p/2$ and $\dmax=q+4=132$, so
\[
c_*(\boldsymbol\pi)=\min\left\{\frac1{\log(4|T|/p)},\frac1{2\log132}\right\}>\frac1{20}>c.
\]
Example~\ref{ex:openai} checks \textup{(H1)--(H3)} for both distributions. Applying Theorem~\ref{thm:partitioned} once to each prescribed root gives, simultaneously with high probability, torsion-freeness, non-elementary hyperbolicity, a finite two-dimensional classifying space, and the root separation used in~\cite[Prop.~5.1]{seed}. Together with the link parity identity of~\cite[(2.1)]{seed}, these are the inputs to the cancellation argument of~\cite[\S6]{seed}, which produces zero divisors in $\mathbb F_2[G_\Phi]$. The zero divisors and root-separation argument originate in that preprint; hyperbolicity of this specific construction is also proved by Mghirbi--Douglas~\cite[Corollary~5.3]{md2026}. The present theorem places this construction within the general link-distribution criterion. Section~\ref{sec:graph-application} gives the cancellation argument in detail.
\end{example}

\subsection{Proof for several parts}

\begin{proof}[Proof of Theorem~\ref{thm:partitioned}]
We use $n=|V_n|$ throughout and check the changes to Sections~\ref{sec:estimates}--\ref{sec:proofs}. All constants may now depend on $\boldsymbol\pi$ and $c$. Write $w_i,P_i$ for the turn probabilities and word weights of $\pi_i$. Set $\wmin=\min_i\wmin(\pi_i)$ and choose common constants $\delta>0$ and $h_0$ in Lemma~\ref{lem:H3forms}(2), taking the smallest of the finitely many decay constants and the largest starting length. Let $i(x)$ be the part containing $x$, and write $V_a(x)=V_a^{(i(x))}$.

First consider the graph estimates. A prescription $(x,s,y)$ between different parts never holds. For a prescription within one part, replace $np_s$ by $|V_s(x)|\ge n\pmin$. In the switching lemma, choose both switched edges in the same part; the switch then stays in $\Omega_n^{\mathrm{pt}}$. In the existence-of-girth and conditioned-prescription arguments, choose the second edge among the edges of the same label in that part. There are $|V_s(x)|$ choices before the same $o(n)$ exclusions. Thus Lemma~\ref{lem:cond-presc} holds with these replacement denominators. The expansion proof uses only the uniform lower bound $|V_s(x)|\ge n\pmin$, so the component-diameter bound also remains valid.

The diagram definitions and their deterministic properties are unchanged. Each trace lies in one part, as does each connected component of the trace graph. Record this part in the pattern, require equivalences to preserve it, and require realizations to send components to their recorded parts. There are at most $K$ components in the trace graph of a diagram with $K$ stars. Recording their parts therefore multiplies the pattern count of Lemma~\ref{lem:count} by at most $J^K$, leaving its polynomial bound unchanged.

In Lemma~\ref{lem:embed}, prescribe a part for each component of $H$. If a non-end vertex $y$ is in part $i$, its image lies in
\[
V_{\inv{a_y}}^{(i)}\cap V_{b_y}^{(i)},
\]
a set of size $|V_{a_y}^{(i)}|w_i(a_y,b_y)$, where $a_y,b_y$ are the entering and leaving labels. The factor from the entering-edge prescription is $|V_{a_y}^{(i)}|^{-1}$, and these factors cancel exactly as before. Each uncancelled factor of this form is at most $(n\pmin)^{-1}$. The fixed-root saving is unchanged. Thus the embedding bound uses the turn matrix of the assigned part at each vertex.

For a piece $\beta$, let $i(\beta)$ be the part of the component containing its chain. Use $P_{i(\beta)}(v_\beta)$ as its weight and
\[
d(v_\beta):=\lfloor-\log P_{i(\beta)}(v_\beta)\rfloor
\]
to record its logarithmic weight. The estimates~\eqref{eq:cost-bounds} and~\eqref{eq:cost-count} hold for each $P_i$ with the common constants chosen above. In the stage bound, send each component of $\Delta_j$ to the part recorded for its containing component of $\Delta$. The turns inside $\beta$ multiply to $P_{i(\beta)}(v_\beta)$, and the initial piece in a constraint component is counted using its own part's word-count estimate~\eqref{eq:cost-count}. The propagation and multiplicity arguments compare letters and pairing constraints and do not require the pieces to have the same turn matrix. The summation uses the same common bounds $\delta$ and $\wmin$. Hence the local inequality follows without further changes. Realized pieces have positive weights because their traces are non-backtracking paths in their assigned graphs.

The deterministic globalization and the proof of Theorem~\ref{intro:construction} now apply with these constants. They use the girth threshold only through the condition $\girth(\Gamma_\Phi)\ge L$ and the inequalities $(c/2)\log n\le L\le c\log n$, valid for all sufficiently large $n$. The same inequalities hold for the thresholds $L_n$ in the statement. This proves the final assertion as well.
\end{proof}

\section{The characteristic-two zero-divisor construction}
\label{sec:graph-application}

We now apply the outlabel criterion to the two-part construction
of~\cite{seed}. The link distributions $\pi_A,\pi_B$ and their collision
estimates are given in Example~\ref{ex:openai}; their common alphabet
has $16516$ signed letters. Example~\ref{ex:partitioned-zero} identifies
the matching law with the partitioned model of
Theorem~\ref{thm:partitioned}. We recall the algebraic argument in full,
because it explains exactly why root separation is sufficient.

The zero-divisor factors and their cancellation are due
to~\cite[Section~6]{seed}. Hyperbolicity for these same groups is also
proved by Mghirbi--Douglas~\cite[Corollary~5.3]{md2026}. The application
below recovers that conclusion from the general outlabel criterion;
we do not claim it as a new existence theorem.

In this section, $n$ denotes the number of vertices in each part,
so the total number is $2n$. Let $A_n,B_n$ have $n$ vertices each,
where $n$ is a multiple of $2(128^2+128+1)^2$, and use the exact
link proportions of Example~\ref{ex:openai}. Write $\tau_{2n}(x)$
for the link assigned to a vertex $x$. The matchings stay within
the two parts.
Condition their union on girth at least
\[
 L_n=\left\lfloor\frac1{100}\log n\right\rfloor,
\]
and write $\mu_{2n}$ for this conditioned measure. The subscript records
the total number of vertices. With $c=1/100$, the difference between
$L_n$ and $c\log(2n)$ is bounded independently of $n$, as allowed by
Theorem~\ref{thm:partitioned}. Fix roots $x_A\in A_n$ and $x_B\in B_n$
before sampling the matchings.

\begin{proposition}[The parity product]\label{g:parity}
Suppose that both roots are separated from every other vertex in
their respective components in the sense of
Theorem~\ref{thm:partitioned}(4). Then
$\F_2[G_\Phi]$ has nonzero elements $\alpha,\beta$ with
$\alpha\beta=0$.
\end{proposition}

\begin{proof}
Let $A',B'$ be the vertex sets of the components containing
$x_A,x_B$, respectively. For $x\in A'$,
let $g_x\in G_\Phi$ be the label of a path from $x_A$ to $x$; define
$h_y$ similarly for $y\in B'$. Closed-path labels are relators, so
these elements do not depend on the chosen paths. Set
\[
 \alpha=\sum_{x\in A'}g_x,
 \qquad
 \beta=\sum_{y\in B'}h_y^{-1}.
\]
Root separation says that $g_x=1$ exactly when $x=x_A$, and $h_y=1$
exactly when $y=x_B$. The identity therefore has coefficient one in
each factor. Thus neither factor vanishes. Distinct vertices away
from a root need not have distinct labels; no full embedding of the
defining graph is required.

For these link distributions,
\begin{equation}\label{g:odd-types}
 |\tau_{2n}(x)\cap\tau_{2n}(y)|\quad\text{is odd whenever }x\in A_n, y\in B_n.
\end{equation}
Indeed, two projective lines over $\F_{128}$ have either one or
$129$ common points. The set of extra letters on the $A$ side has
size zero or two; on the $B$ side it is empty or consists of all
three extra letters. Their intersection has size zero or two.
Adding these contributions proves~\eqref{g:odd-types}.

Form the finite graph with vertex set $A'\times B'$ whose oriented
edges are simultaneous steps in $A'$ and $B'$ with the same signed
label. A label $t$ can be used at $(x,y)$ precisely when
$t\in\tau_{2n}(x)\cap\tau_{2n}(y)$. The immersion property makes the step
unique, and reversing it replaces $t$ by $\bar t$. Thus the degree
of $(x,y)$, with loops counted twice, is
$|\tau_{2n}(x)\cap\tau_{2n}(y)|$, which is odd. The handshaking lemma implies
that every connected component of this graph has an even number
of vertices.

If a simultaneous step takes $(x,y)$ to $(x',y')$ with label $t$,
then $g_{x'}=g_xt$ and $h_{y'}=h_yt$, whence
\[
 g_{x'}h_{y'}^{-1}=(g_xt)(h_yt)^{-1}=g_xh_y^{-1}.
\]
The product $g_xh_y^{-1}$ is therefore constant on each connected
component. In the expansion
\[
 \alpha\beta=\sum_{(x,y)\in A'\times B'}g_xh_y^{-1},
\]
each component contributes an even number of copies of a single
group element. Every such contribution is zero over $\F_2$, proving
$\alpha\beta=0$.
\end{proof}

\begin{corollary}[Generic hyperbolic examples in the two-part model]
\label{g:genericity-theorem}
Under the conditioned law $\mu_{2n}$, with probability tending
to one as $n$ tends to infinity through the stated multiples,
$G_\Phi$ is torsion-free and non-elementary hyperbolic,
$X_\Phi$ is a finite two-dimensional classifying space, and
$\F_2[G_\Phi]$ has nonzero zero divisors. The last assertion also
holds over every field of characteristic two.
\end{corollary}

\begin{proof}
The link distributions satisfy (H1)--(H3), and the girth parameter
is allowed by Example~\ref{ex:partitioned-zero}. Apply
Theorem~\ref{thm:partitioned} once with root $x_A$ and once with root
$x_B$. The intersection of these two events still has probability
tending to one, by the union bound. It supplies hyperbolicity,
the finite classifying space, torsion-freeness, and both instances
of root separation. Proposition~\ref{g:parity} supplies the zero
divisors. If $k$ has characteristic two, the coefficient inclusion
$\F_2\hookrightarrow k$ induces an injective map
$\F_2[G_\Phi]\hookrightarrow k[G_\Phi]$, preserving the two nonzero
factors and their product.
\end{proof}

The probability here is on matching data, with exact link counts
and a girth condition. It is not a uniform measure on isomorphism
classes, nor the independent-relator density law. There is no claim
that examples at different sizes are pairwise nonisomorphic.

Theorem~\ref{thm:partitioned}(2) also gives a quantitative conclusion about the marking: with probability tending to one, $F(S)\to G_\Phi$ is injective on the ball of radius
$\lfloor\varepsilon_1L_n/2\rfloor$. If two distinct elements of
that ball had the same image,
cyclic reduction of their quotient would produce a nonempty
word of length at most $\varepsilon_1L_n$ that is trivial in $G_\Phi$.
Consequently these marked groups converge locally, in probability,
to the free group as $n$ grows.

The parity calculation is specific to characteristic two. Mghirbi and Douglas obtain the stronger failure of direct finiteness in every positive characteristic; for odd primes they use signed coefficients and different vertex types \cite[Main Theorem and \S\S6--7]{md2026}.

\section{What changes in the standard density model?}\label{sec:comparison}

In the outlabel model, closed paths supply relators, so shared paths impose correlations among them. We record two restrictions on zero divisors when relators are instead sampled as independent words. They show why the genericity conclusion of Corollary~\ref{g:genericity-theorem} requires a separate analysis in the standard density model.

Fix a set $X$ of $r\ge2$ generators, let $Q=2r-1$, and let $\mathcal W_\ell$ be the set of cyclically reduced words of length $\ell$ in $X^{\pm1}$, counted as literal strings. For $0<d<1/2$, put $N_\ell=\lfloor Q^{d\ell}\rfloor$ and choose $w_1,\ldots,w_{N_\ell}$ independently and uniformly from $\mathcal W_\ell$, with replacement. Write $G_{r,d,\ell}=\langle X\mid w_1,\ldots,w_{N_\ell}\rangle$. Throughout this section, limits are taken as $\ell\to\infty$, with $r$ and $d$ fixed. We use the identity
\[
 |\mathcal W_\ell|=Q^\ell+r+(r-1)(-1)^\ell.
\]
This is the trace of the $\ell$th power of the letter-transition matrix, whose eigenvalues are $Q$, $1$ with multiplicity $r$, and $-1$ with multiplicity $r-1$.

\subsection{Division rings below density one sixth}

\begin{theorem}[Division rings at low density]\label{c:lowdensity}
For every fixed $r\ge2$ and $0<d<1/6$, with probability tending to one
$G_{r,d,\ell}$ is torsion-free, hyperbolic, and virtually compact special.
On this event, $k[G_{r,d,\ell}]$ embeds in a division ring for every
division ring $k$. In particular, for every field $k$,
\[
 \Prob\bigl[k[G_{r,d,\ell}]\text{ has nonzero zero divisors}\bigr]
 \longrightarrow0,
\]
and the probability that $k[G_{r,d,\ell}]$ is directly finite tends to one.
\end{theorem}

\begin{proof}
Random groups at these densities are torsion-free and hyperbolic with
probability tending to one~\cite[Section~2]{ollivier}.
Ollivier--Wise~\cite{ow} prove that they act freely and cocompactly on
CAT(0) cube complexes with probability tending to one. Agol's
theorem~\cite[Theorem~1.1]{agol} then makes them virtually compact special.
Fisher--S\'anchez-Peralta~\cite[Theorem~1.3(1)]{fsp} embed the group algebra
of every torsion-free virtually compact special group over every division
ring into a division ring. Such an embedding excludes nonzero zero
divisors and implies direct finiteness: if $ab=1$, then the same identity
holds in the containing division ring, where it forces $ba=1$.

These conclusions hold for independent sampling as well as for sampling
a uniform subset of literal cyclically reduced words. Under independent
sampling, the probability of a repeated word is at most
\[
 \frac{\binom{N_\ell}{2}}{|\mathcal W_\ell|}
 =O\bigl(Q^{(2d-1)\ell}\bigr)=o(1).
\]
Conditioned on no repetition, the underlying set is a uniform
$N_\ell$-element subset of $\mathcal W_\ell$. This comparison transfers
events of probability tending to one between the two sampling conventions.

\end{proof}

\begin{corollary}[A restriction on cubical actions]\label{c:cubical}
Let $G$ be torsion-free. Suppose that, for some field $k$, the group
algebra $k[G]$ has nonzero zero divisors or is not directly finite.
Then $G$ is not virtually compact special. If $G$ is also hyperbolic,
it admits no proper
cocompact action by cubical automorphisms on a CAT(0) cube complex.
\end{corollary}

\begin{proof}
The division-ring embedding of~\cite{fsp} excludes both algebraic
phenomena for torsion-free virtually compact special groups. A proper
cocompact cubical action of a hyperbolic group would imply virtual
compact specialness by~\cite[Theorem~1.1]{agol}.
\end{proof}

This is a consequence of the cited results, not a new obstruction.
Mghirbi--Douglas apply the same argument to their positive-characteristic
examples in~\cite[Corollary~7.13 and Section~10]{md2026}.
The inference from a cubical action uses cocompactness; no claim is
made here about arbitrary proper cubical actions.

\subsection{Fixed support sizes at every subcritical density}

The next obstruction concerns support sizes, not the lengths of words
representing support elements. Its input is the first-order theory of
random groups. For every fixed existential sentence $\sigma$ in the
language of groups, every fixed $r\ge2$, and $0<d<1/2$, Massalha proves
\cite[Theorem~10.3]{massalha} that
\begin{equation}\label{c:firstorder}
 \lim_{\ell\to\infty}\Prob[G_{r,d,\ell}\models\sigma]
 =\begin{cases}
 1,&F_r\models\sigma,\\
 0,&F_r\not\models\sigma.
 \end{cases}
\end{equation}
Kharlampovich--Miasnikov--Sklinos~\cite[Theorem~2]{kms} prove the stronger
statement for arbitrary first-order sentences. The existential case
already suffices below.

The model of independent cyclically reduced words is explicit
in~\cite[Definition~2.1]{massalha}, with
$\lfloor|\mathcal W_\ell|^d\rfloor$ samples. The counting estimate in
its Lemma~3.9 uses the upper bound $|\mathcal W_\ell|^{dt}$ for choices
of $t$ sample indices. The same estimates apply to
$N_\ell=\lfloor Q^{d\ell}\rfloor$, since $Q^\ell\le|\mathcal W_\ell|$.

\begin{proposition}[Fixed-support obstruction over a fixed field]
\label{c:support}
Let $k$ be any fixed field. For every fixed $M\ge1$, $r\ge2$, and
$0<d<1/2$,
\[
 \Prob\!\left[
 \begin{gathered}
  \exists\,0\ne\alpha,\beta\in k[G_{r,d,\ell}],\qquad\alpha\beta=0,\\
  |\supp\alpha|\le M,\qquad |\supp\beta|\le M
 \end{gathered}\right]\longrightarrow0.
\]
There is no restriction on the word lengths of the support elements.
\end{proposition}

\begin{proof}
Fix positive integers $a,b$. Write $I_{a,b}=\{1,\ldots,a\}\times
\{1,\ldots,b\}$. Call a partition $\mathcal P$ of this finite set
\emph{$k$-admissible} if there are coefficients
$\lambda_1,\ldots,\lambda_a,\mu_1,\ldots,\mu_b\in k\setminus\{0\}$
such that
\begin{equation}\label{c:admissible}
 \sum_{(i,j)\in C}\lambda_i\mu_j=0
 \qquad\text{for every block }C\in\mathcal P.
\end{equation}
Let $\mathfrak P_{a,b}(k)$ be the set of these partitions. It is finite,
even when $k$ is infinite. For each of its members, fix one coefficient
tuple satisfying~\eqref{c:admissible}. No effective procedure for finding
the admissible partitions is needed.

For a partition $\mathcal P$, let
\[
 E_{\mathcal P}(x,y):=
 \bigwedge_{C\in\mathcal P}\;
 \bigwedge_{(i,j),(i',j')\in C}x_i y_j=x_{i'}y_{j'}.
\]
Consider the following existential sentence, containing only group
variables:
\begin{equation}\label{c:sentence}
\begin{split}
 \sigma_{a,b,k}:=\exists x_1\cdots x_a y_1\cdots y_b\;\Bigg[
 &\bigwedge_{i<i'}x_i\ne x_{i'}\ \wedge\
 \bigwedge_{j<j'}y_j\ne y_{j'}\\
 &\wedge\bigvee_{\mathcal P\in\mathfrak P_{a,b}(k)}
 E_{\mathcal P}(x,y)\Bigg].
\end{split}
\end{equation}
An empty disjunction is understood to be false. The subscript $k$
records the finite choice of allowed partitions; it does not add field
variables or field operations to the language of groups.

We claim that $G\models\sigma_{a,b,k}$ if and only if $k[G]$ has
nonzero factors with support sizes exactly $a$ and $b$ whose product is zero.
Suppose first that a disjunct $E_{\mathcal P}$ holds. Use the coefficient
tuple chosen for $\mathcal P$ and set
\[
 \alpha=\sum_{i=1}^a\lambda_i x_i,
 \qquad \beta=\sum_{j=1}^b\mu_j y_j.
\]
The distinctness conditions and nonzero coefficients give the required
support sizes. Within each block of $\mathcal P$, all products $x_i y_j$
are equal and their total coefficient is zero. Different blocks may
give the same group element; merging their zero coefficient sums still
gives zero. Thus $\alpha\beta=0$. This explains why no inequalities
between different blocks are required.

Conversely, write any such zero product with its distinct support
elements and nonzero coefficients. Partition $I_{a,b}$ according to
equality of the group elements $x_i y_j$. The vanishing of each
coefficient of $\alpha\beta$ says exactly that this partition is
$k$-admissible. Its equalities give a true disjunct
of~\eqref{c:sentence}. This proves the claim.

The free group $F_r$ is bi-orderable. In an order invariant under left
and right multiplication, the largest product of elements of two
finite nonempty supports is attained only by the pair of largest
support elements. Its coefficient is therefore a product of two
nonzero field elements and cannot vanish. Hence $k[F_r]$ has no nonzero zero divisors,
so $F_r\not\models\sigma_{a,b,k}$. Equation~\eqref{c:firstorder} gives
$\Prob[G_{r,d,\ell}\models\sigma_{a,b,k}]\to0$. Take the finite union
over $1\le a,b\le M$.
\end{proof}

Over $\F_p$, the admissible partitions can also be found by enumerating
the finitely many tuples of nonzero coefficients. Over $\F_2$,
all nonzero coefficients are $1$, so every product value must occur
an even number of times. Pairing those occurrences recovers the usual
characteristic-two formulation. The partition argument above is what
allows the statement to cover arbitrary fixed fields.

There is a corresponding fixed-support conclusion for direct
finiteness. If $a,b\in k[G]$ satisfy $ab=1$ but $ba\ne1$, put
$z=1-ba$. Then $z\ne0$, $a\ne0$, and $az=0$. If both $a$ and $b$
have support size at most $M$, then $z$ has support size at most
$M^2+1$. Proposition~\ref{c:support}, applied with this larger fixed
bound, therefore makes the probability of such a direct-finiteness
witness tend to zero as well.

The unrestricted zero-divisor event is a union over all support bounds
$M$. One cannot interchange this infinite union with the limit without
an additional uniform estimate.
Thus Proposition~\ref{c:support} does not determine the probability of
the unrestricted event for $1/6\le d<1/2$. Any generic existence statement
would have to allow witnesses that escape every fixed common bound on
the two support sizes.

\section*{Use of artificial intelligence}
\addcontentsline{toc}{section}{Use of artificial intelligence}
ChatGPT (OpenAI) was used in developing and reviewing arguments,
checking calculations and dependencies, locating and comparing
references, and preparing the exposition and \LaTeX{} source. The
outlabel material also incorporates a draft developed with Claude
(Anthropic). This assistance is not independent mathematical or
formal verification. External results and their authors are credited
at the points where they are used.

\end{document}